\documentclass[12pt]{amsart}

\usepackage{graphicx} 
\usepackage[margin=1in]{geometry}

\usepackage{verbatim}
\usepackage{amssymb}
\usepackage{amscd}
\newtheorem{lemma}{Lemma}[section]
\newtheorem{proposition}[lemma]{Proposition}
\newtheorem{theorem}[lemma]{Theorem}

\theoremstyle{definition}
\newtheorem{defn}[lemma]{Definition}

\theoremstyle{remark}
\newtheorem{remark}[lemma]{Remark}
\newtheorem{example}[lemma]{Example}

\newcommand{\Z}{\ensuremath{{\mathbb Z}}}
\newcommand{\B}{\ensuremath{\mathbb{B}}}
\newcommand{\N}{\ensuremath{{\mathbb N}}}

\newcommand{\R}{\ensuremath{\mathbb R}}
\newcommand{\C}{\ensuremath{\mathbb C}}

\renewcommand{\Pr}{\ensuremath{\mathbb P}}

\newcommand{\A}{\ensuremath{\mathcal{A}}}
\newcommand{\F}{\ensuremath{\mathcal{F}}}
\newcommand{\K}{\ensuremath{\mathcal{K}}}
\newcommand{\T}{\ensuremath{\mathcal{T}}}

\author{Nadia Alluhaibi}
\address{N. Alluhaibi,  Science and Arts College, 
Rabigh Campus, 
King Abdulaziz University,
Jeddah, Saudi Arabia }
\email{nallehaibi@kau.edu.sa}

\author{Tatyana Barron}
\address{T. Barron, Department of Mathematics, 
University of Western Ontario,
London, Ontario N6A 5B7, Canada }

\email{tatyana.barron@uwo.ca}

\thanks{Research of T. Barron is supported in part 
by the Natural Sciences and Engineering Research Council of Canada}

\title{On vector-valued automorphic forms on bounded symmetric domains}

\begin{document}

\sloppy

\maketitle

\noindent {\bf Abstract.} We prove a spanning result for vector-valued Poincar\'e series on a
bounded symmetric domain. 
We associate a sequence of  holomorphic automorphic forms to a submanifold of the domain.
When the domain is the unit ball in $\C ^n$, 
we provide estimates for the norms of these automorphic forms and we find 
asymptotics of the norms  (as the weight goes to infinity)
for a class of totally real submanifolds. We give an example of a CR submanifold
of the ball, for which the norms of the associated automorphic forms have 
a different asymptotic behaviour. 

\


\


\section{Introduction}

Let $D$ be a bounded symmetric domain (in $\C^n$, for $n\ge 1$).
Suppose $\Gamma$ is a cocompact discrete 
subgroup of $Aut(D)$. Let $k$ be a positive integer. 
A holomorphic function $f:D\to\C$ is called a {\it holomorphic 
automorphic form of weight $k$ for} $\Gamma$ if $f(\gamma z)J(\gamma , z)^k=f(z)$ for all $\gamma\in \Gamma$, $z\in D$. 
Here $J(\gamma ,z)$ denotes the determinant of the Jacobi matrix of $\gamma$ at $z$. 
Let $m$ be a positive integer and let $\rho:\Gamma\to GL(m,\C )$ be a unitary representation of $\Gamma$. 
A holomorphic function $F:D\to\C^m$ is called a {\it holomorphic $\C ^m$-valued  
automorphic form of weight $k$} (for the pair ($\Gamma$, $\rho$)) 
if $F(\gamma z)=\rho(\gamma)F(z)J(\gamma , z)^{-k}$ for all $\gamma\in \Gamma$, $z\in D$. 
Holomorphic automorphic forms correspond to holomorphic sections of ${\mathcal{L}}^{\otimes k}$, where ${\mathcal{L}}$ 
is the canonical bundle on 
$M=\Gamma\backslash D$, and $\C ^m$-valued holomorphic automorphic forms correspond 
to holomorphic sections of $E_{\rho}\otimes {\mathcal{L}}^{\otimes k}$, 
where $E_{\rho}\to M$ is the flat vector bundle defined by $\rho$. 

The theory of automorphic forms is a vast subject that has strong interaction with 
many areas of mathematics, including 
representation theory, number theory, semiclassical analysis and quantization.  
One connection 
between automorphic forms and quantization is as follows: $M$ is a K\"ahler manifold,
${\mathcal{L}}$ is a quantum line bundle, $\frac{1}{k}$ is interpreted as $\hbar$ (the Planck constant),
and the space of holomorphic automorphic forms, with the Petersson inner product, is isomorphic to the Hilbert space 
$H^0(M,{\mathcal{L}}^{\otimes k})$ 
used in quantization \cite{barron:12}. Berezin-Toeplitz quantization or K\"ahler quantization is usually studied for $\C$-valued observables.  
Extending the theory to $\C^m$-valued functions on $D$ is a non-trivial task which is interesting
from the mathematical 
point of view and physically meaningful (see e.g. work by S.T. Ali and M. Englis \cite{ali:07} on domains 
in $\C^n$).

There are many different kinds of automorphic/modular forms, and generalizations. Vector-valued automorphic forms are 
ubiquitous and go back to classical works of Borel, Selberg and others (see, for example \cite{borel:65}). 
Applications include work of R. Borcherds on singular Howe correspondence, 
work of S. Kudla on arithmetic cycles, physics-related work by T. Gannon, G. Mason and others. The space of Jacobi 
forms is isomorphic to a space of vector-valued modular forms. Various kinds of 
vector-valued forms for $G=SU(n,1)$ (i.e. when $D$ is the $n$-dimensional complex hyperbolic space, or, equivalently, the open unit ball in $\C^n$ with the complex 
hyperbolic metric: $D=\B^n\simeq SU(n,1)/S(U(n)\times U(1))$) 
have been studied, in particular, in recent papers by E. Freitag, G. van der Geer and others \cite{clery:13}, \cite{freitag:14}, 
in work by Kato including \cite{kato:84}, and in work by Kojima including \cite{kojima:97}. 
It is well known that modular forms appear in generating functions for arithmetic or algebraic 
objects in Calabi-Yau varieties. 
It would be interesting to see if Picard modular forms could play a similar role. 

Poincar\'e series is a standard and powerful
tool that is used in automorphic forms, spectral theory, complex analysis,
Teichm\"uller theory, algebraic geometry, and other areas.
In \cite{foth:07, foth:08} T. Barron (Foth) studied automorphic forms 
for compact smooth $M=\Gamma \backslash D$  (in \cite{foth:08} $D=\B^n$),  
constructing explicitly the automorphic form $f_p$ ($p\in D$) with the property $(g,f_p)=g(p)$ 
for any other holomorphic automorphic form $g$. Here $(.,.)$ denotes the Petersson inner product. 
Such $f_p$ is constructed via Poincar\'e series and is related to the (weighted) Bergman kernel 
and to the concept of a coherent state. Choosing a $q$-form on a $q$-dimensional submanifold, and 
integrating $f_p$, 
one can get automorphic forms associated to submanifolds of $D$. In this paper we extend this 
framework to vector-valued holomorphic automorphic forms. 
In Section \ref{secspan}
we prove that sufficiently many of these vector-valued Poincar\'e series
span the space of vector-valued 
holomorphic automorphic forms  - Theorem \ref{thspan}. 

There are somewhat different, but closely related results in literature: 
it is known that $\C$-valued Poincar\'e series of polynomials in $z_1$,...,$z_n$
span the space of holomorphic automorphic forms
on a bounded symmetric domain (for sufficiently large weights)
\cite{bell:67, foth:07, wu:15}. 
In David Bell's thesis \cite{bell:67} it is stated that similar results also hold
for vector-valued automorphic forms on classical domains and it is explained how the proofs for $\C$-valued case
can be extended to the vector-valued case. 

To give general context to Section \ref{secsubman}, we observe that   
associating an automorphic form or, more generally, a section of a vector bundle,
to a submanifold 
of a K\"ahler manifold is an idea that is used in many contexts.
In particular, relative Poincar\'e series 
can be associated to closed geodesics on a hyperbolic Riemann surface
\cite{katok:85, katok:87}. S. Katok and  
T. Foth (Barron)  generalized this construction
from compact Riemann surfaces of genus $g\ge 1$ 
to ball quotients in \cite{foth:01, foth:04} (where they addressed the spanning problem),
and more recently  
T.B. addressed the non-vanishing question in \cite{barron:18}.  
In \cite{kudla:79, tong:83} the submanifold is a closed geodesic or, more generally, a totally 
geodesic submanifold. To mention a somewhat different kind of such technique, 
there is a way to associate a section of a line bundle 
to a Bohr-Sommerfeld Lagrangian submanifold, 
which is used in semiclassical analysis and symplectic geometry (see, in particular,
\cite{borthwick:95, burns:10, deber:06, gorod:01, jeffrey:92, 
paoletti:08}). In 
\cite{guillemin:15, ioos:18} sections of vector bundles are associated to  isotropic submanifolds.

Here, in Section \ref{secsubman},
we take advantage of the fact that the K\"ahler manifold is $M=\Gamma\backslash D$, 
the holomorphic sections of the vector bundle on $M$ correspond to holomorphic vector-valued
functions on $D$, 
and we associate an automorphic form to a submanifold of a fundamental domain of $\Gamma$ in $D$, 
and not to a submanifold of $M$ (see Remark \ref{bohrs}). 
In the case when the domain is the unit ball in $\C^n$, we provide asymptotic
(as the weight goes to infinity) 
statements about the inner products  - Theorems \ref{thxy}, \ref{thxx}.
In particular, 
in Theorem \ref{thxy}{\it (ii)}
we show that if two submanifolds are at a positive distance from each other, then
the inner product of the associated automorphic forms decreases rapidly as the weight goes
to infinity. In Theorem \ref{thxx}{\it (ii)} we show, in particular,
that for a class of totally real submanifolds,  as $k\to\infty$, 
the square of the norm of $\Theta_X^{(j;k)}$
grows as a positive constant
times $k^{n-\frac{q}{2}}$, where
$X$ is such a submanifold,  $\Theta_X^{(j;k)}$ is one of the $m$ $\C^m$-valued Poincar\'e series
associated to $X$ ($j\in \{ 1,2,...,m\}$), and  
$q$ is the real dimension of $X$.  
We work out several examples. In Example \ref{3ball} we estimate the asymptotics
of the norms for a $3$-dimensional submanifold which is CR and not totally real, and we show
that the leading term in the asymptotics the square of the norm of the associated Poincar\'e
series is \underline{not}  
$const \cdot k^{n-\frac{3}{2}}$.

This paper contains results from the Ph.D. thesis of N.A. \cite{alluhaibi:17}
written under the supervision of T.B.

{\bf Acknowledgments.} We are thankful to A. Dhillon, Y. Karshon,  M. Pinsonnault, E. Schippers,
A. Uribe and N. Yui for related discussions. We acknowledge the referee's efforts.

\section{Preliminaries}
\label{secprelim}

Let $D=G/K\subset \C^n$, for $n\ge 1$, be a bounded symmetric domain
($G=Aut(D)$ is a real semisimple
Lie group that acts transitively on $D$, $K$ a maximal compact subgroup of $G$).
Denote by $z_1$,...,$z_n$ the complex coordinates. Also denote 
$z_j=x_j+iy_j$ ($x_j,y_j\in\R$, $1\le j\le n$) and  
denote the Euclidean volume form by 
$$
dV_e=dx_1\wedge dy_1\wedge ... \wedge dx_n\wedge dy_n=
\Bigl (\frac{i}{2}\Bigr ) ^n dz_1\wedge d\bar{z}_1\wedge ... \wedge dz_n\wedge d\bar{z}_n.
$$   
Let $\K (.,.)$ be the Bergman kernel for $D$. It has the reproducing property: 
$$
f(z)=\int\limits_D f(w) \K (z,w)dV_e (w),
$$
$z\in D$, for all functions $f$ that are holomorphic on $D$ and such that $\int\limits_D |f(z)|^2dV_e(z)<\infty$. Also 
$\K (z,w)=\overline{\K (w,z)}$ for $z,w\in D$, and 
\begin{equation}
\label{bktransf}  
J(\gamma,z) \overline{J(\gamma,w)}\K (\gamma z,\gamma w)=\K (z,w)
\end{equation}
for $z,w\in D$, $\gamma \in G$, 
where $J(\gamma,z)$ is the complex Jacobian of the transformation $D\to D$ at $z$ defined by $\gamma$.
The $(1,1)$-form 
 $\omega =i\partial\bar{\partial}\log K(z,z)$ is a $G$-invariant K\"ahler form on $D$.

Let $k\in\N$ be a positive integer. It will be usually assumed that $k$ is sufficiently large. 
The volume form $dV(z)=\K (z,z) dV_e(z)$ is $G$-invariant. 
The reproducing kernel for the Hilbert space 
of holomorphic functions on $D$ satisfying $\int\limits_D |f(z)|^2\K (z,z)^{-k}dV(z)<\infty$ 
is $c(D,k)\K (z,w)^k$, where $c(D,k)$ is a constant (note: the proof of this fact in
\cite{selberg:57} uses the assumption 
that $Aut(D)$ acts transitively on $D$). 
The reproducing property is: for any such function $f$ 
\begin{equation}
  \label{reprod}
f(z)=c(D,k)\int\limits_D f(w) \K (z,w)^k\K (w,w)^{-k}dV (w),
\end{equation}
$z\in D$. 
The value of the constant is determined by 
(1.6)\cite{selberg:57}: 
\begin{equation}
\label{constcdk}
c(D,k)\int\limits_D\frac{\K (z,w)^k \K (w,z)^k}{\K (w,w)^k}dV(w)=\K (z,z)^k
\end{equation}
for any $z\in D$.

Let $\Gamma$ be a discrete subgroup of $G$ such that the quotient 
$M=\Gamma \backslash D=\Gamma\backslash G/K$ is smooth and compact. 
Let $m$ be a positive integer.     
Let $\rho : \Gamma\to GL(m,\C )$ be a unitary representation of $\Gamma$. 
\begin{defn} \cite{baily:73}
A function $f:D\to\C$ is called a {\it (holomorphic) $\Gamma$-automorphic 
form of weight $k$} if $f$ is holomorphic and 
\begin{equation}
  \label{autfs}
f(\gamma z)J(\gamma ,z)^k=f(z) \ \forall 
\gamma  \in \Gamma, z\in D .
\end{equation}
\end{defn} 
\begin{defn} \cite{selberg:57} 
\label{vvafdef} 
A {\it vector-valued automorphic form of weight $k$ for $(\rho, \Gamma )$} 
is $F=\begin{pmatrix} F_1 \\ ... \\ F_m \end{pmatrix}$, where $F_j: D\to\C$, $j=1,...,m$, are holomorphic functions, 
and 
\begin{equation}
\label{autfvv}
J(\gamma,z)^kF(\gamma z)=\rho(\gamma)F(z) \  \forall 
\gamma   \in \Gamma, z\in D . 
\end{equation}
\end{defn}
Denote the space of holomorphic $\Gamma$-automorphic forms of weight $k$ on $D$
by $\A(\Gamma,k)$. Denote the space of $\C ^m$-valued
holomorphic $(\rho,\Gamma)$-automorphic forms of weight $k$ on $D$
by $\A(\Gamma,m,k,\rho)$.
\begin{remark} 
In a more general case when $M$ is of finite volume and not compact the definitions should include an appropriate 
condition at the cusps.
The condition "$M$ is smooth" can be relaxed to allow $\Gamma$ such as, for example,  
$SL(2,\Z)\subset SL(2,\R)\simeq SU(1,1)$ 
or $SU(2,1)\cap SL(3,\Z [i])$. 
\end{remark} 

Define the inner product on the space $\A(\Gamma,m,k,\rho)$ as follows: 
\begin{equation}
  \label{inprod}
( F,G) = \int\limits_{\Gamma \backslash D} F(z)^T\overline{ G(z)} \K (z,z)^{-k} dV
\end{equation}
for $F,G\in \A(\Gamma,m,k,\rho)$. 
This is well-defined because the function $F(z)^T\overline{ G(z)} \K (z,z)^{-k}$ is $\Gamma$-invariant 
(note: for that it is essential that $\rho$ is unitary).  

Define the inner product on the space $\A(\Gamma,k)$ by  
$$
( f,g) = \int\limits_{\Gamma \backslash D} f(z)\overline{ g(z)} \K (z,z)^{-k} dV
$$
for $f,g\in \A(\Gamma,k)$.

Denote by $K_M$ the canonical bundle on $M$ and by $K_D$ the canonical bundle on $D$. 
\begin{remark} We have isomorphisms of Hilbert spaces: 
  $\A(\Gamma,k)\cong H^0(M,K_M^{\otimes k})$,  
  $\A(\Gamma,m,k,\rho)\cong H^0(M,E_{\rho}\otimes K_M^{\otimes k})$.
  In particular, a holomorphic function $f$ on $D$ satisfies (\ref{autfs}) if and only if
  $f(z)(dz_1\wedge ...\wedge dz_n)^{\otimes k}$ is a $\Gamma$-invariant holomorphic section of
  $K_D^{\otimes k}$ (and thus descends to a holomorphic section of $K_M^{\otimes k}$).  
\end{remark} 
\begin{remark}
  There are irreducible unitary representations of the fundamental group
  of a compact Riemann surface of genus $\ge 2$ for each $m\in\N$ (Proposition 2.1 
  \cite{narasimhan:64}). The proof in \cite{narasimhan:64} provides explicit examples of such
  representations. 
\end{remark}

\section{Poincar\'e series and a spanning result}
\label{secspan}

Let $D$ be a bounded symmetric domain, and let $\Gamma$ be a discrete subgroup of $Aut(D)$
such that the quotient $M=\Gamma\backslash D$ is smooth and compact. 
Let $k$ and $m$ be positive integers, and let $\rho$ be an $m$-dimensional unitary representation of $\Gamma$.
This is the setting for this section. 

For an integrable holomorphic function $F:D\to\C^m$ we define, formally, the Poincar\'e series
of weight $k$ 
$$
\Theta_F(z)=\sum_{\gamma\in\Gamma} \rho(\gamma ^{-1})F(\gamma z)J(\gamma , z)^k
$$
(here we omit $k$ from notation and write simply $\Theta_F$). 
If the series converges uniformly on compact sets in $D$, 
then $\Theta_F\in \A(\Gamma,m,k,\rho)$. Indeed,
since the convergence is uniform on compact sets,
it follows that $\Theta_F$ is holomorphic. To verify (\ref{autfvv}), we observe:   
for $g\in \Gamma$, $z\in D$ 
$$
\Theta_F(g z)=\sum_{\gamma\in\Gamma} \rho(\gamma ^{-1})F(\gamma gz)J(\gamma , gz)^k= 
\sum_{\gamma\in\Gamma} \rho(g(\gamma g) ^{-1})F(\gamma gz)\frac{J(\gamma g, z)^k}{J(g, z)^k}= 
\rho(g) J(g, z)^{-k}\Theta_F(z).
$$ 
Choose $p\in D$.  
In \cite{foth:07} the $\C$-valued Poincar\'e series 
$$
\theta_p(z)=\sum_{\gamma\in\Gamma}\Bigl ( \K (\gamma z,p)J(\gamma ,z)\Bigr ) ^k \in \A(\Gamma,k)
$$ 
(convergent absolutely and uniformly on compact sets for sufficiently large $k$, and  
having  the property    
$(f,\theta_p)={\mbox{const}}(D,k)f(p)$
for any $f\in \A(\Gamma,k)$) were considered, and it was shown 
that such Poincar\'e series for an appropriate number of points in general position 
form a basis in $\A(\Gamma,k)$. Note that the property $(f,\theta_p)={\mbox{const}}(D,k)f(p)$ reflects 
the fact that the Bergman kernel for $K_M^{\otimes k}$ is the Poincar\'e series of the Bergman kernel for 
$K_D^{\otimes k}$ (Theorem 2 \cite{ma:15} or Theorem 1 \cite{lu:15}).

Let us now generalize the construction from  \cite{foth:07} 
by associating to a point $p\in D$ $m$ vector-valued Poincar\'e series 
\begin{equation}
\label{vvps}
\Theta_p^{(j;k)}(z)=c(D,k) \sum_{\gamma\in\Gamma}\rho(\gamma^{-1}) \T _p (\gamma z)  J(\gamma,z)^k, \ j=1,...,m 
\end{equation}
where $\T _p (z)=\begin{pmatrix} (\T_p) _1 (z) \\ ... \\ (\T _p) _m (z) \end{pmatrix}$, 
$(\T _p)_j(z)=\K ( z,p)^k$ and $(\T _p)_l (z)=0$ for $l\ne j$  
(i.e. $\T _p (\gamma z)$ is the vector-function whose components, except for the $j$-th one, are zero, 
and $(\T _p)_j(\gamma z)=  \K ( \gamma z,p)^k$). 

We shall also use the notation $\hat{\Theta}^{(j;k)}(z,p)$ for the function 
$$
\hat{\Theta}^{(j;k)}:D\times D\to \C^m 
$$
$$
(z,p)\to \Theta_p^{(j;k)}(z).
$$  
\begin{lemma}
\label{propconvscalar}
Let $p\in D$.  
For $k\ge 2$ 
the series $\sum\limits_{\gamma\in\Gamma}(\K (\gamma z,p)J(\gamma,z))^k$ converges absolutely and 
uniformly on compact sets of $D$.  
\end{lemma}
The proof is in the Appendix. 
\begin{proposition}
  \label{propconv}
  Let $j\in \{ 1,...,m\}$ and let $p\in D$. 
Suppose $k$ is  sufficiently large. 

\noindent (i) The series (\ref{vvps}) converges absolutely and uniformly on compact sets. 

\noindent (ii) For each $H\in \A(\Gamma,m,k,\rho)$
$$
( H,\Theta_p^{(j;k)})=H_j(p). 
$$
\end{proposition}
\begin{theorem}
\label{thspan}
For sufficiently large $k$,  
for sufficiently many points $p_1,...,p_d$ in general position, the $\C$-linear span of 
$\{ \Theta_{p_l}^{(j;k)} | 1\le l\le d; \  1\le j\le m\}$ is 
$\A(\Gamma,m,k,\rho)$. 
\end{theorem}
\noindent {\bf Proof of Proposition \ref{propconv}.} Proof of {\it (i)}. 
For $1\le l\le m$
$$
\Bigl | \Bigl ( \rho(\gamma^{-1}) \T _p (\gamma z) J(\gamma,z)^k
\Bigr )_l \Bigr |\le 
\sqrt{\Bigl (\rho(\gamma^{-1}) \T _p (\gamma z) J(\gamma,z)^k 
\Bigr )^T \overline{\rho(\gamma^{-1}) \T _p (\gamma z) J(\gamma,z)^k} }=
$$
$$
|\K (\gamma z,p)J(\gamma,z)|^k.
$$
The statement now follows from Lemma \ref{propconvscalar}.

Proof of {\it (ii)}. Let $\F$ be a Dirichlet fundamental domain for $\Gamma$ 
(or a canonical fundamental domain \cite{selberg:68}).  
Denote $w=\gamma z$ for $\gamma\in\Gamma$, $z\in \F$. By (\ref{autfvv}) 
 $H(z)^T=H(w)^T(\rho (\gamma) ^{-1})^T J(\gamma,z)^k$.
Using (\ref{bktransf}), (\ref{reprod}), (\ref{inprod}), (\ref{vvps}), we get:  
$$
( H,\Theta_p^{(j;k)}) =c(D,k) \int\limits_{\F} \begin{pmatrix} H_1(z) & ... & H_m(z) \end{pmatrix} 
\sum_{\gamma\in\Gamma}\overline{\rho(\gamma^{-1})} \overline{\T _p(\gamma z) } \overline{J(\gamma,z)^k}
\K (z,z)^{-k}dV(z)=
$$
$$
c(D,k)\sum_{\gamma\in\Gamma}\int\limits_{\gamma \F} H(w)^T \overline{\T _p(w)}   K(w,w)^{-k}dV(w)= 
$$
$$
c(D,k)\sum_{\gamma\in\Gamma}\int\limits_{\gamma \F} H_j(w)\overline{\K (w,p)}^k K(w,w)^{-k}dV(w)=
$$
$$
c(D,k)\int\limits_{D} H_j(w) \K (p,w)^k K(w,w)^{-k}dV(w)=
H_j(p).  
$$ 
$\Box$ 

\noindent {\bf Proof of Theorem \ref{thspan}}.
Let $k\in\N$. 
The holomorphic vector bundle $W=E_{\rho}\otimes K_M^{\otimes k}$ is positive.
Using the notations similar to those in Chapters 2, 3 \cite{kobayashi:87},  
denote by $P(W)$ the fibre bundle over $M$ whose fibre at $x$ is $\Pr (W_x)$ 
(i.e. $P(W)=(W-\{ \text{zero section} \})/\C^*$), denote by $\pi: P(W)\to M$ the projection, 
and by $L(W)$ the tautological line bundle over $P(W)$ 
(i.e. the subbundle of $\pi^*W$ with the fiber $L(W)_\xi$ 
at $\xi\in P(W)$ being the complex line in $W_{\pi(\xi)}$ represented by $\xi$).  
Also denote by $L(W^*)$ the tautological line bundle over $P(W^*)=(W^*-\{ \text{zero section} \})/\C^*$
and by $\hat{\pi}: P(W^*)\to M$ the projection. 
We note that a section $s$ of $W$  produces a section $\tilde{s}$ of $(L(W^*))^*\to P(W^*)$.
Specifically, $\tilde{s}=h\circ s\circ \hat{\pi}$, 
where $h$ is the holomorphic surjection $\hat{\pi}^*W\to (L(W^*))^*\simeq L(W)$ given, fiberwise, 
by the quotient map $W_x\to W_x/ \ker f$ over $(x,[f])\in P(W^*)$, where $x\in M$, $f\in W_x^*$, $f\ne 0$. 

Suppose $k$ is large enough, so that the line bundle $(L(W^*))^*\to P(W^*)$ is very ample. 
Let $d=\dim H^0(P(W^*),(L(W^*))^*)$ and let $\tilde{p}_1$,...,$\tilde{p}_d$ be points 
in $P(W^*)$ in general position (i.e. such that their images under the projective embedding given 
by  $(L(W^*))^*$ are not on the same hyperplane in $\Pr(H^0(P(W^*),(L(W^*))^*))^*$).
Such $d$ points exist because the linear system is base point free.  
Select a Dirichlet fundamental domain $\F$ for $\Gamma$ and for each $j\in \{ 1,...,d\}$ let 
$p_j$ be the point in $\F$ that corresponds to $\hat{\pi}(\tilde{p}_j)$. 

Now, to prove the statement of the theorem, 
suppose $H\in \A(\Gamma,m,k,\rho)$ is not in the linear span of
$\{ \Theta_{p_l}^{(j;k)}\} _{\substack{1\le l\le d \\ 1\le j\le m}}$. 
Then $H$ is in the orthogonal complement of this subspace and therefore  
$$
( H,\Theta_{p_l}^{(j;k)})=0
$$
for all $l\in \{ 1,...,d\}$ and all $j\in \{ 1,...,m\}$.
By Proposition \ref{propconv}{\it (ii)} $H(p_1)=...=H(p_d)=0$.  
Let $s$ be the section of $W$ corresponding to $H$. 
This section vanishes at $p_1$,...,$p_d$. 
Therefore $\tilde{s}(\tilde{p_1})=...=\tilde{s}(\tilde{p_d})=0$. Since $\tilde{p_1}$,...,$\tilde{p_d}$ 
are in general position, we conclude that $\tilde{s}\equiv 0$. It follows that 
 $s=0$. Hence $H=0$.  
$\Box$

\section{Automorphic forms and submanifolds}
\label{secsubman}

Let $D$ be a bounded symmetric domain, and let $\Gamma$ be a discrete subgroup of $Aut(D)$
such that the quotient $M=\Gamma\backslash D$ is smooth and compact. 
Let $k$ and $m$ be positive integers, and let $\rho$ be an $m$-dimensional unitary representation of $\Gamma$.

\noindent Let $\Lambda$ be a  $q$-dimensional submanifold of $D$ ($q\ge 1$) such that 
$\Lambda \subset \bar{B}(z_0,r_0) \subset D$, where $\bar{B}(z_0,r_0)$ is the closed ball centered at $z_0$ 
of radius $r_0$ with respect to the Euclidean metric, 
for some $z_0\in D$, $r_0>0$.    
Let $\nu$ be a nonzero volume form (a real $q$-form)  on $\Lambda$
such that $\int\limits_{\Lambda}\nu >0$.  
Set  
\begin{equation}
\label{setup1}
\Theta_{\Lambda}^{(j;k)}(z)=\int\limits_{\Lambda} \hat{\Theta}^{(j;k)}(z,p)\K (p,p)^{-\frac{k}{2}}\nu(p) 
\end{equation}
for $j\in \{1,...,m\}$.
By a standard differentiation under the integral sign argument $\Theta_{\Lambda}^{(j;k)}$ is holomorphic. 
Moreover, 
$\Theta_{\Lambda}^{(j;k)}\in \A(\Gamma,m,k,\rho)$ and 
\begin{equation}
\label{setup2}
(H,\Theta_{\Lambda}^{(j;k)})=\int\limits_{\Lambda} H_j(z)\K (z,z)^{-\frac{k}{2}}\nu(z)
\end{equation}
for any $H\in \A(\Gamma,m,k,\rho)$.
\begin{remark}
  The statement analogous to 
  Proposition \ref{propconv}{\it (ii)}, but  written for
  the corresponding sections of $E_{\rho}\otimes K_M^{\otimes k}$, 
  would mean that the section of  $E_{\rho}\otimes K_M^{\otimes k}$, corresponding  to
  $\Theta_{p}^{(j;k)}$,  
  is  the $j$-th row of the Bergman kernel for this vector bundle,
  where the Bergman kernel is written as an $m\times m$  matrix.
  The general idea of ``integrating the Bergman kernel over a submanifold $\Lambda$''
  was used in \cite{borthwick:95} (to obtain sections of powers of a line bundle,
  with $\Lambda$ being a Bohr-Sommerfeld Lagrangian submanifold of a compact K\"ahler
  manifold) and it is used in a recent preprint \cite{ioos:18}
  (to obtain sections of certain vector bundles, with $\Lambda$ being
  an isotropic  Bohr-Sommerfeld submanifold of a compact symplectic manifold). 
\end{remark}  
In this section,
the domain $D$ will be the unit ball $\B^n\subset \C^n$ ($n\ge 1$), with its Bergman metric.
 Recall that
$SU(n,1)=\{ A\in SL(n+1,\C ) \ | \ 
A^T\sigma \bar{A}=\sigma \}$, where 
$\sigma =\begin{pmatrix} 
1_{n\times n} & 0 \cr 0 & -1\end{pmatrix}$. The ball 
is a bounded realization of the Hermitian symmetric space $SU(n,1)/S(U(n)\times U(1))$ 
(note that for $n=1$ $D$ is the unit disc: 
$D\cong SU(1,1)/U(1)\cong SL(2,\R)/SO(2)$). The group $SU(n,1)$ acts on 
$\B^n$ is by fractional-linear transformations:  
for $\gamma =(a_{jk})\in SU(n,1)$ the corresponding automorphism  
$\B^n\to \B^n$ is 

\noindent 
$z=(z_1,...,z_n)\mapsto \Bigl ( \dfrac{a_{11}z_1+...+a_{1n}z_n+a_{1,n+1}}{a_{n+1,1}z_1+...
+a_{n+1,n}z_n+a_{n+1,n+1}} , ..., 
\dfrac{a_{n1}z_1+...+a_{nn}z_n+a_{n,n+1}}{a_{n+1,1}z_1+...
+a_{n+1,n}z_n+a_{n+1,n+1}} \Bigr )$.

The complex Jacobian is 
$J(\gamma,z)=1/(a_{n+1,1}z_1+...+a_{n+1,n}z_n+a_{n+1,n+1})^{n+1}$.
\begin{remark}
  Each element of the center of $SU(n,1)$ acts as the identity map on $\B^n$, and $Aut(\B^n)$ is
  isomorphic to $PU(n,1)$. We will represent automorphisms of the ball by matrices from $SU(n,1)$,
  and we will use the same letter to denote the matrix and the corresponding automorphism. 
\end{remark}
We will denote by $0$ the point 
$(0,...,0)\in \B^n$. 
Also, for $z,w\in \B^n$ denote 
$$
\langle z,w\rangle=z_1\bar{w}_1+...+z_n\bar{w}_n-1.
$$
The $SU(n,1)$-invariant K\"ahler form on $\B^n$ is, up to a positive constant factor, 
$$
i\partial\bar{\partial}\log (-\langle z,z\rangle)=
\frac{i}{\langle z,z\rangle ^2}\Bigl [
  (\sum_{j=1}^n\bar{z}_jdz_j)\wedge (\sum_{l=1}^nz_ld\bar{z}_l)-\langle z,z\rangle
  \sum_{r=1}^ndz_r\wedge d\bar{z}_r\Bigr ].
$$
Denote by $\tau (z,w)$ the distance between $z$ and $w$ with respect 
to the complex hyperbolic metric. Note that 
\begin{equation}
  \label{coshdistf}
\cosh ^2\frac{\tau(z,w)}{2}=\frac{\langle z,w\rangle \langle w,z\rangle}{\langle z,z\rangle \langle w,w\rangle}
\end{equation}
(see e.g. \cite{goldman:99} 3.1.7). 
It is a standard fact (see e.g. \cite{rudin:80} or \cite{range:86}) that for the ball 
\begin{equation}
  \label{bergkball}
\K (z,w)=\frac{n!}{\pi^n}(-\langle z,w\rangle)^{-(n+1)}.
\end{equation}
\begin{lemma}
\label{lemconstcdk}
For $D=\B^n$ the constant $c(D,k)$ from Section \ref{secprelim} is $c(\B^n ,k)=\binom{(n+1)(k-1)+n}{n}$. 
\end{lemma}
This follows from Theorem 2.2 \cite{zhu:05} with $\alpha=(n+1)(k-1)$ (the constant $c(D,k)$ comes out to be 
$c_{\alpha}$ given by (2.2)\cite{zhu:05}). This also can be verified in another way, by a direct calculation (see  
the Appendix). 
\begin{remark}
\label{remconstball}
Applying the Stirling formula $N!\sim (\frac{N}{e})^N\sqrt{2\pi N}\Bigl (1+O(\frac{1}{N})\Bigr )$ as $N\to\infty$ \cite{debru:61},  
we get:  $c(\B^n ,k)\sim \frac{(n+1)^n}{n!}k^n\Bigl (1+O(\frac{1}{k})\Bigr )$ as $k\to\infty$.
\end{remark}
Let $\Gamma$ be a discrete subgroup of $SU(n,1)$ such that the quotient 
$M=\Gamma \backslash \B^n$ is smooth and compact.
Let $k$, $m$ be positive integers, and let $\rho$ be an $m$-dimensional unitary representation of $\Gamma$.  
Denote by $\pi:\B^n\to M$ the covering map. Let $\F$ be a Dirichlet fundamental domain for $\Gamma$ \cite{parker:09}.  
Suppose $X$ and $Y$ are submanifolds of $\B^n$ of dimensions $q_X>0$ 
and $q_Y>0$ respectively, such that $X=\pi^{-1}(X')\cap \F$, $X\cong X'$, and  
$Y=\pi^{-1}(Y')\cap \F$, $Y\cong Y'$, where $X'$ and $Y'$ are submanifolds of $M$, and
$\cong$ stands for diffeomorphism.   
Let $\nu_X$ be a nonzero volume form on $X$ (a real $q_X$-form) such that $\int_X \nu_X>0$  
and let $\nu_Y$ is a nonzero volume form on $Y$ (a real $q_Y$-form) such that $\int_Y \nu_Y>0$.  
Denote $\tilde{X}=\Gamma X$, $\tilde{Y}=\Gamma Y$. 
Define the $q_X$-form $\nu_{\tilde{X}}$ on ${\tilde{X}}$ by 
$\nu_{\tilde{X}}\Bigr |_{\gamma^{-1}(X)}=\gamma^*\nu_X$ for each $\gamma\in\Gamma$. 
Define $\nu_{\tilde{Y}}$ the same way. Note that $\nu_{\tilde{X}}$, $\nu_{\tilde{Y}}$ 
are $\Gamma$-invariant. 
Assume 
$\int\limits_{\tilde{X}}|K(z,w)|^2\frac{\nu_{\tilde{X}}(w)}{K(w,w)}<\infty$ for all $z\in \F$, 
$\int\limits_{\tilde{Y}}|K(z,w)|^2\frac{\nu_{\tilde{Y}}(w)}{K(w,w)}<\infty$ for all $z\in \F$ 
(the last condition is satisfied, for example, when $Y$ is a small ball and $\nu_{\tilde{Y}}=dV\Bigr | _{\tilde{Y}}$, because
 $K(.,w)$ is square-integrable on $\B^n$).  
 
For two subsets $A$, $B$ of $\B^n$ we will denote
$$
{\mathrm{dist}}(A,B)=\inf \{ \tau(z,w) \ | z\in A, w\in B\}.
$$
Since $\tau$ is $\Gamma$-invariant, the same notation can be used for subsets $A$, $B$ of $M$. 
\begin{remark}
Recall that if $(a_k)$, $(b_k)$ are two sequences of complex numbers, then notation $a_k\sim b_k$ as $k\to\infty$ means 
$\lim\limits_{k\to\infty}\frac{a_k}{b_k}=1$.  
\end{remark}
\begin{theorem}  
\label{thxy} Let $r,j\in \{ 1,...,m\}$. 

\noindent (i) Suppose ${\mathrm{dist}}(\tilde{X}-X,Y)>0$  and $r\ne j$.  Then 
for any $l\in\N$ there is a constant $C=C(l;n,X,Y,\Gamma , \nu_X, \nu_Y)$ 
such that, as $k\to\infty$ 
$$
|(\Theta_X^{(r;k)},\Theta_Y^{(j;k)})|\le \frac{C}{k^l}. 
$$ 
(ii) Suppose ${\mathrm{dist}}(\tilde{X},Y)>0$. 
Then 
for any $l\in\N$ there is a constant $C=C(l;n,X,Y,\Gamma , \nu_X, \nu_Y)$ 
such that, as $k\to\infty$ 
$$
|(\Theta_X^{(r;k)},\Theta_Y^{(j;k)})|\le \frac{C}{k^l}. 
$$   
\end{theorem} 
\begin{remark}
  \label{remdist}
  If ${\mathrm{dist}}(X,\partial \F)>0$ or ${\mathrm{dist}}(Y,\partial \F)>0$,
  then ${\mathrm{dist}}(\tilde{X}-X,Y)>0$. 
\end{remark}
\begin{remark}
  If  ${\mathrm{dist}}(\tilde{X},Y)>0$, then ${\mathrm{dist}}(X',Y')>0$. 
\end{remark}  
\begin{theorem}  
\label{thxx} Suppose $Y\subset X$, ${\mathrm{dist}} (X,\partial\F)>0$, and $j\in \{1,...,m\}$.  

\noindent (i) If $q_X\le n$, then  
$$
(\Theta_X^{(j;k)},\Theta_Y^{(j;k)})\le {\mbox{const}}(n,X,Y,\nu_X,\nu_Y)k^{n-\frac{q_X}{2}}
$$ 
as $k\to\infty$. 

\noindent (ii) If $X\subset \{ z\in\B^n | y_1=...=y_n=0\}$, then  
$$
(\Theta_X^{(j;k)},\Theta_Y^{(j;k)})\sim C(n,X,Y,\nu_X,\nu_Y)k^{n-\frac{q_X}{2}} 
$$ 
as $k\to\infty$, where $C(n,X,Y,\nu_X,\nu_Y)$ is a positive constant. 
\end{theorem}
\begin{remark}
  \label{bohrs}
  Because of the assumptions in Theorem \ref{thxx}{\it (ii)}, in this part
  of the theorem the submanifolds $X$ and $Y$ are isotropic submanifolds
  of $\B^n$. In Theorem \ref{thxx}{\it (i)} and in Theorem \ref{thxy}
  the submanifolds are not necessarily isotropic. Note that $X$ and $Y$ are submanifolds of $\B^n$, the
  universal cover of $M$, and not of $M$. We do not require $X$ and $Y$ to satisfy a Bohr-Sommerfeld
  condition. In the usual procedure of associating 
  a section of a line bundle, $L$, to a Lagrangian or isotropic 
  submanifold, $\Lambda$, of  a K\"ahler manifold, the Bohr-Sommerfeld condition ensures the existence
  of a covariant constant nonvanishing section, $\varphi$, of $L^*\Bigr |_{\Lambda}$.
  Having such $\varphi$, from a holomorphic section $s$ of $L^{\otimes k}$, one obtains  
  a function on $\Lambda$, $\varphi^{\otimes k}(s)$, which then can be integrated over $\Lambda$. This   
  provides a linear functional on the space of holomorphic sections of $L^{\otimes k}$.   
  We do not need such $\varphi$, since in (\ref{setup2}) we are already integrating a function. 
\end{remark}
\noindent {\bf Proof of Theorem \ref{thxy}.}
Using (\ref{vvps}), (\ref{setup1}), (\ref{setup2}), we get: 
$$
|(\Theta_X^{(r;k)},\Theta_Y^{(j;k)})|=|\int\limits_{Y}(\Theta_X^{(r;k)}(z))_j \K ( z,z)^{-\frac{k}{2}}\nu_{Y}(z)|=
$$
$$
|\int\limits_{Y}\int\limits_{X}(\hat{\Theta}^{(r;k)}(z,\zeta))_j \K ( \zeta,\zeta)^{-\frac{k}{2}} \nu_{X}(\zeta)
\K ( z,z)^{-\frac{k}{2}} \nu_{Y}(z)|=
$$
$$
c(\B^n,k)|\int\limits_{Y}\int\limits_{X}\sum_{\gamma\in\Gamma}\rho(\gamma^{-1})_{jr}\K (\gamma z,\zeta)^kJ(\gamma,z)^k 
\K ( \zeta ,\zeta)^{-\frac{k}{2}} \nu_{X}(\zeta)
\K (z,z)^{-\frac{k}{2}} \nu_{Y}(z) |\le
$$  
$$
c(\B^n,k)\int\limits_{Y}\int\limits_{X}\sum_{\gamma\in\Gamma}|\rho(\gamma^{-1})_{jr}||\K (\gamma z,\zeta)J(\gamma,z)|^k 
\K ( \zeta ,\zeta)^{-\frac{k}{2}} \nu_{X}(\zeta)
\K (z,z)^{-\frac{k}{2}} \nu_{Y}(z) .
$$  
Setting $\zeta=\gamma w$, we get, using (\ref{bktransf}): 
$$
|(\Theta_X^{(r;k)},\Theta_Y^{(j;k)})|\le c(\B^n,k)\int\limits_{Y}\sum_{\gamma\in\Gamma}
\int\limits_{\gamma ^{-1} X}
|\rho(\gamma^{-1})_{jr}||\K (z,w)|^k \K (w,w)^{-\frac{k}{2}} \nu_{\tilde{X}}(w)
\K (z,z)^{-\frac{k}{2}} \nu_{Y}(z).
$$
Since $\rho(\gamma^{-1})$ is a unitary matrix, we have: $|\rho(\gamma^{-1})_{jr}|\le 1$. 
Using (\ref{coshdistf}), (\ref{bergkball}), we get: 
$$
|(\Theta_X^{(r;k)},\Theta_Y^{(j;k)})|\le
c(\B^n,k)\int\limits_{Y}\sum_{\gamma\in\Gamma}
\int\limits_{\gamma ^{-1} X}
|\K (z,w)|^k \K (w,w)^{-\frac{k}{2}} \nu_{\tilde{X}}(w)
\K (z,z)^{-\frac{k}{2}} \nu_{Y}(z)=
$$
$$
c(\B^n,k)\int\limits_{Y}
\int\limits_{\tilde{X}}
|\K (z,w)|^k \K (w,w)^{-\frac{k}{2}} \nu_{\tilde{X}}(w)
\K (z,z)^{-\frac{k}{2}} \nu_{Y}(z)=
$$
$$
c(\B^n ,k)\int\limits_{Y}\int\limits_{\tilde{X}}
|\K (z,w)|^k \K (w,w)^{-\frac{k}{2}+1} \frac{\nu_{\tilde{X}}(w)}{K(w,w)}
\K (z,z)^{-\frac{k}{2}+1}  \frac{\nu_{Y}(z)}{K(z,z)}=
$$ 
$$
c(\B^n ,k)  \int\limits_{Y}\int\limits_{\tilde{X}}
\Bigl ( \frac{\langle z,z\rangle\langle w,w\rangle }{\langle z,w\rangle  \langle w,z\rangle} 
\Bigr )^{ (n+1)(\frac{k}{2}-1)  }
|K(z,w)|^2    \frac{\nu_{\tilde{X}}(w)}{K(w,w)} \frac{\nu_{Y}(z)}{K(z,z)}=  
$$
$$
c(\B^n ,k)  \int\limits_{Y}\int\limits_{\tilde{X}}
\Bigl ( \cosh \frac{\tau(z,w)}{2}\Bigr ) ^{-(n+1)(k-2)}
|K(z,w)|^2   \frac{\nu_{\tilde{X}}(w)}{K(w,w)} \frac{\nu_{Y}(z)}{K(z,z)}   \le 
$$
$$
c(\B^n ,k)  \Bigl ( \frac{1}{\cosh [\frac{1}{2}{\mathrm{dist}}(\tilde{X},Y)] }\Bigr )^ {(n+1)(k-2)}
\int\limits_{Y}\int\limits_{\tilde{X}} |K(z,w)|^2   \frac{\nu_{\tilde{X}}(w)}{K(w,w)} \frac{\nu_{Y}(z)}{K(z,z)}  .
$$
For $r\ne j$, since $\rho({\mathrm{id}})_{jr}=0$, the argument above can be modified:  
$$
|(\Theta_X^{(r;k)},\Theta_Y^{(j;k)})|\le c(\B^n,k)\int\limits_{Y}\sum_{\gamma\in\Gamma, \gamma\ne id}
\int\limits_{\gamma ^{-1} X}
|\K (z,w)|^k \K (w,w)^{-\frac{k}{2}} \nu_{\tilde{X}}(w)
\K (z,z)^{-\frac{k}{2}} \nu_{Y}(z)=
$$
$$
c(\B^n,k)\int\limits_{Y}
\int\limits_{\tilde{X}-X}
|\K (z,w)|^k \K (w,w)^{-\frac{k}{2}} \nu_{\tilde{X}}(w)
\K (z,z)^{-\frac{k}{2}} \nu_{Y}(z)=
$$
$$
c(\B^n ,k)\int\limits_{Y}\int\limits_{\tilde{X}-X}
|\K (z,w)|^k \K (w,w)^{-\frac{k}{2}+1} \frac{\nu_{\tilde{X}}(w)}{K(w,w)}
\K (z,z)^{-\frac{k}{2}+1}  \frac{\nu_{Y}(z)}{K(z,z)}=
$$ 
$$
c(\B^n ,k)  \int\limits_{Y}\int\limits_{\tilde{X}-X}
\Bigl ( \frac{\langle z,z\rangle\langle w,w\rangle }{\langle z,w\rangle  \langle w,z\rangle} 
\Bigr )^{ (n+1)(\frac{k}{2}-1)  }
|K(z,w)|^2    \frac{\nu_{\tilde{X}}(w)}{K(w,w)} \frac{\nu_{Y}(z)}{K(z,z)}=  
$$
$$
c(\B^n ,k)  \int\limits_{Y}\int\limits_{\tilde{X}-X}
\Bigl ( \cosh \frac{\tau(z,w)}{2}\Bigr ) ^{-(n+1)(k-2)}
|K(z,w)|^2   \frac{\nu_{\tilde{X}}(w)}{K(w,w)} \frac{\nu_{Y}(z)}{K(z,z)}   \le 
$$
$$
c(\B^n ,k)  \Bigl ( \frac{1}{\cosh [\frac{1}{2}{\mathrm{dist}}(\tilde{X}-X,Y)] }\Bigr )^ {(n+1)(k-2)}
\int\limits_{Y}\int\limits_{\tilde{X}-X} |K(z,w)|^2   \frac{\nu_{\tilde{X}}(w)}{K(w,w)} \frac{\nu_{Y}(z)}{K(z,z)}  .
$$
Now,
$$
0<\int\limits_{Y}\int\limits_{\tilde{X}-X} |K(z,w)|^2 \frac{\nu_{\tilde{X}}(w)}{K(w,w)} \frac{\nu_{Y}(z)}{K(z,z)}\le
$$
$$
\int\limits_{Y}\int\limits_{\tilde{X}} |K(z,w)|^2 \frac{\nu_{\tilde{X}}(w)}{K(w,w)} \frac{\nu_{Y}(z)}{K(z,z)}  ={\mbox{const}} (\tilde{X},Y,\nu_{\tilde{X}},\nu_Y)
$$  
Since $\cosh \frac{\varepsilon}{2}>1$ for $\varepsilon>0$, and with Remark \ref{remconstball}, 
the statements follow.  $\Box$

\noindent {\bf Proof of Theorem \ref{thxx}.} 
Using (\ref{vvps}), (\ref{setup1}), (\ref{setup2}), we get: 
$$
(\Theta_X^{(j;k)},\Theta_Y^{(j;k)})= \int\limits_{Y}(\Theta_X^{(j;k)}(z))_j \K ( z,z)^{-\frac{k}{2}}\nu_{Y}(z)=
$$
$$
\int\limits_{Y}\int\limits_{X}(\hat{\Theta}^{(j;k)}(z,\zeta))_j \K ( \zeta,\zeta)^{-\frac{k}{2}} \nu_{X}(\zeta)
\K ( z,z)^{-\frac{k}{2}} \nu_{Y}(z)=
$$
$$
c(\B^n ,k)\int\limits_{Y}\int\limits_{X}\sum_{\gamma\in\Gamma}\rho(\gamma^{-1})_{jj}(\K (\gamma z,\zeta)J(\gamma,z))^k 
\K ( \zeta ,\zeta)^{-\frac{k}{2}} \nu_{X}(\zeta)
\K (z,z)^{-\frac{k}{2}} \nu_{Y}(z) = 
$$
$$
I_1+I_2,
$$
where $I_1$ is the term with $\gamma={\mbox{id}}$ and $I_2$ is the rest. Thus,  
$$
I_1=c(\B^n ,k)\int\limits_{Y}\int\limits_{X}\K (z,\zeta)^k 
\K ( \zeta ,\zeta)^{-\frac{k}{2}} \nu_{X}(\zeta)
\K (z,z)^{-\frac{k}{2}} \nu_{Y}(z)
$$
and 
$$
I_2=c(\B^n ,k)\int\limits_{Y}\int\limits_{X}
\sum_{\gamma\in\Gamma, \gamma\ne {\mbox{id}}}\rho(\gamma^{-1})_{jj}(\K (\gamma z,\zeta)J(\gamma,z))^k 
\K ( \zeta ,\zeta)^{-\frac{k}{2}} \nu_{X}(\zeta)
\K (z,z)^{-\frac{k}{2}} \nu_{Y}(z). 
$$
Since $\rho(\gamma^{-1})$ is a unitary matrix, $|\rho(\gamma^{-1})_{jj}|\le 1$. 
Setting $\zeta=\gamma w$ and using (\ref{bktransf}), (\ref{coshdistf}), (\ref{bergkball}), 
we get: 
$$
|I_2|\le 
c(\B^n ,k)\int\limits_{Y}\int\limits_{X}\sum_{\gamma\in\Gamma, \gamma\ne {\mbox{id}} }|\K (\gamma z,\zeta)J(\gamma,z)|^k 
\K ( \zeta ,\zeta)^{-\frac{k}{2}} \nu_{X}(\zeta)
\K (z,z)^{-\frac{k}{2}} \nu_{Y}(z) =
$$
$$ 
 c(\B^n ,k)\int\limits_{Y}\sum_{\gamma\in\Gamma , \gamma\ne {\mbox{id}} } \ \int\limits_{\gamma ^{-1} X}
|\K (z,w)|^k \K (w,w)^{-\frac{k}{2}} \nu_{\tilde{X}}(w)
\K (z,z)^{-\frac{k}{2}} \nu_{Y}(z)=
$$
$$
c(\B^n ,k)\int\limits_{Y}\int\limits_{\tilde{X}-X}
|\K (z,w)|^k \K (w,w)^{-\frac{k}{2}} \nu_{\tilde{X}}(w)
\K (z,z)^{-\frac{k}{2}}  \nu_{Y}(z)=
$$
$$
c(\B^n ,k)
\int\limits_{Y}\int\limits_{\tilde{X}-X}
\Bigl ( \frac{\langle z,z\rangle\langle w,w\rangle }{\langle z,w\rangle  \langle w,z\rangle} 
\Bigr )^{ (n+1)(\frac{k}{2}-1)  }|K(z,w)|^2  \frac{\nu_{\tilde{X}}(w)}{K(w,w)} \frac{\nu_{Y}(z)}{K(z,z)}=
$$
$$
c(\B^n ,k)
\int\limits_{Y}\int\limits_{\tilde{X}-X}
\Bigl ( \cosh \frac{\tau(z,w)}{2}\Bigr )  ^{ -(n+1)(k-2) }
|K(z,w)|^2  \frac{\nu_{\tilde{X}}(w)}{K(w,w)} \frac{\nu_{Y}(z)}{K(z,z)}\le
$$
$$
c(\B ^n,k)
\Bigl ( \cosh [\frac{1}{2}{\mathrm{dist}}(Y,\tilde{X}-X)]\Bigr )  ^{ -(n+1)(k-2) }
\int\limits_{Y}\int\limits_{\tilde{X}-X}
|K(z,w)|^2  \frac{\nu_{\tilde{X}}(w)}{K(w,w)} \frac{\nu_{Y}(z)}{K(z,z)}.
$$
Because of Remark \ref{remdist}, $\cosh [\frac{1}{2}{\mathrm{dist}}(Y,\tilde{X}-X)]>1$. 
Also  
$$
\int\limits_{Y}\int\limits_{\tilde{X}-X}
|K(z,w)|^2  \frac{\nu_{\tilde{X}}(w)}{K(w,w)} \frac{\nu_{Y}(z)}{K(z,z)}\le 
\int\limits_{Y}\int\limits_{\tilde{X} }
|K(z,w)|^2  \frac{\nu_{\tilde{X}}(w)}{K(w,w)} \frac{\nu_{Y}(z)}{K(z,z)}<\infty.
$$
Using Remark \ref{remconstball}, 
we conclude that $I_2$ has the property: for any $l\in\N$ there is a constant
$C=C(l;n,X,Y,\Gamma , \nu_X,\nu_Y)$ 
such that 
$$
|I_2|\le \frac{C}{k^l}  
$$ 
as $k\to\infty$.

Now we consider
$$
I_1=
c(\B^n ,k)  \int_Y\int_X
 \frac{(\langle z,z\rangle\langle \zeta,\zeta\rangle)^{ \frac{(n+1)k}{2}  } }
{(-\langle z,\zeta\rangle )^{(n+1)k} } 
\nu_{X}(\zeta)\nu_{Y}(z).
$$
We use Fubini's theorem to switch  
to the integral over $Y\times X$ with respect to the product measure, then choose and fix 
a sufficiently small $\delta>0$, and split $I_1$ into two parts: 
$I_1^{(1)}$, where the integration is over the part of  $Y\times X$ where $\tau(z,\zeta)\le \delta$ 
and  $I_1^{(2)}$, where the integration is over the part of  $Y\times X$ where $\tau(z,\zeta)> \delta$.  
Using (\ref{coshdistf}), we get: 
$$
I_1^{(2)}=
c(\B^n ,k)  
\iint\limits_{\substack{Y\times X \\ \tau(z,\zeta)> \delta} }
 \frac{(\langle z,z\rangle\langle \zeta,\zeta\rangle)^{ \frac{(n+1)k}{2}  } }
{(-\langle z,\zeta\rangle )^{(n+1)k} } 
\nu_{Y}(z)\nu_{X}(\zeta), 
$$
$$
|I_1^{(2)}|\le 
c(\B^n ,k)  
\iint\limits_{\substack{Y\times X \\ \tau(z,\zeta)> \delta} }
\Bigl ( \frac{\langle z,z\rangle\langle \zeta,\zeta\rangle }{\langle z,\zeta\rangle  \langle \zeta,z\rangle} 
\Bigr )^{ \frac{(n+1)k}{2}  }\nu_{Y}(z)\nu_{X}(\zeta) = 
$$
$$
c(\B^n ,k)  
\iint\limits_{\substack{Y\times X \\ \tau(z,\zeta)> \delta} }
\Bigl (\cosh \frac{\tau (z,\zeta)}{2}\Bigr )^{-(n+1)k}\nu_{Y}(z)\nu_{X}(\zeta)\le 
$$
$$
c(\B^n ,k)  \frac{1}{ \Bigl ( \cosh \frac{\delta}{2}\Bigr )^{(n+1)k} }
\iint\limits_{\substack{Y\times X \\ \tau(z,\zeta)> \delta} }\nu_{Y}(z)\nu_{X}(\zeta), 
$$
therefore by Remark \ref{remconstball}, and since $\cosh \frac{\delta}{2}>1$,
it follows  that $I_1^{(2)}$ has the property: 
for any $l\in\N$ there is a constant $C=C(l;n,X,Y,\delta , \nu_X, \nu_Y)$ 
such that 
$$
|I_1^{(2)}|\le \frac{C}{k^l} 
$$ 
as $k\to\infty$.

It remains to investigate the term 
$$
I_1^{(1)}=
c(\B^n ,k)  
\iint\limits_{\substack{Y\times X \\ \tau(z,\zeta)\le \delta} }
\frac{(\langle z,z\rangle\langle \zeta,\zeta\rangle) ^{ \frac{(n+1)k}{2}  }}
{(-\langle z,\zeta\rangle )  
^{ (n+1)k } }\nu_{Y}(z)\nu_{X}(\zeta).  
$$ 
To proceed with the proof of  $(i)$, we observe:  
$$
|I_1^{(1)}|\le 
c(\B^n ,k)  
\iint\limits_{\substack{Y\times X \\ \tau(z,\zeta)\le \delta} }
\Bigl ( \frac{\langle z,z\rangle\langle \zeta,\zeta\rangle }{\langle z,\zeta\rangle  \langle \zeta,z\rangle} 
\Bigr )^{ \frac{(n+1)k}{2}  }\nu_{Y}(z)\nu_{X}(\zeta).  
$$
For the proof of $(ii)$: if $z\in Y$ and $\zeta\in X$, then $\langle z,\zeta\rangle =\langle \zeta ,z\rangle$ and 
$$
I_1^{(1)}=
c(\B^n ,k)  
\iint\limits_{\substack{Y\times X \\ \tau(z,\zeta)\le \delta} }
\Bigl ( \frac{\langle z,z\rangle\langle \zeta,\zeta\rangle }{\langle z,\zeta\rangle  \langle \zeta,z\rangle} 
\Bigr )^{ \frac{(n+1)k}{2}  }\nu_{Y}(z)\nu_{X}(\zeta). 
$$
Thus, to finish the proof of the theorem we need to treat the integral
$$
c(\B^n ,k)  
\iint\limits_{\substack{Y\times X \\ \tau(z,\zeta)\le \delta} }
\Bigl ( \frac{\langle z,z\rangle\langle \zeta,\zeta\rangle }{\langle z,\zeta\rangle  \langle \zeta,z\rangle} 
\Bigr )^{ \frac{(n+1)k}{2}  }\nu_{Y}(z)\nu_{X}(\zeta). 
$$
Using (\ref{coshdistf}), we get: 
$$
\iint\limits_{\substack{Y\times X \\ \tau(z,\zeta)\le \delta} }
\Bigl ( \frac{\langle z,z\rangle\langle \zeta,\zeta\rangle }{\langle z,\zeta\rangle  \langle \zeta,z\rangle} 
\Bigr )^{ \frac{(n+1)k}{2}  }\nu_{Y}(z)\nu_{X}(\zeta)=  
\iint\limits_{\substack{Y\times X\\ \tau(z,\zeta)\le \delta} }
\Bigl (\cosh \frac{\tau (z,\zeta)}{2}\Bigr )^{-(n+1)k}\nu_{Y}(z)\nu_{X}(\zeta)=
$$
$$
\iint\limits_{\substack{Y\times X \\ \tau(z,\zeta)\le \delta} }
e^{ -(n+1)k\ln \cosh \frac{\tau (z,\zeta)}{2} }\nu_{Y}(z)\nu_{X}(\zeta)=
\int\limits_{Y}\int\limits_{ \{ \zeta\in X|\tau(z,\zeta)\le \delta \}  }
e^{ -(n+1)k\ln \cosh \frac{\tau (z,\zeta)}{2} }\nu_{X}(\zeta)\nu_{Y}(z).
$$
Let $A_z\in SU(n,1)$, $z \in Y$, be a smooth family of automorphisms 
$\B^n\to \B^n$ such that $A_z z=0$. Denote $\hat{X}=\cup_{z\in Y} A_z(X)$.  
Let $\{ U_j\}$ be a finite cover of $\hat{X}$ 
by open subsets of $\B^n$ with smooth boundary, let $t_1^{(j)}$,...,$t_{q_X}^{(j)}$ be local 
coordinates on $U_j\cap \hat{X}$, and let $\{ \psi^{(j)}\}$ be a partition of unity subordinate to the cover $\{ U_j\}$.

For a fixed $z\in Y$ consider the integral   
$$
\int\limits_{ \{ \zeta\in X|\tau(z,\zeta)\le \delta \}  }
e^{ -(n+1)k\ln \cosh \frac{\tau (z,\zeta)}{2} }\nu_{X}(\zeta)=
\int\limits_{ \{ w\in A_z(X)|\tau(0,w)\le \delta \}  }
e^{ -(n+1)k\ln \cosh \frac{\tau (w,0)}{2} }[(A_z^{-1})^*\nu_{X}](w), 
$$
where $w=A_z \zeta$. Note: 
$\tau(0,w)=\tau(A_zz,A_z\zeta )=\tau(z,\zeta)$.
We have:  $(A_z^{-1})^*\nu_{X}\Bigr | _{U_j}=f^{(j)}(t)  dt_1^{(j)}\wedge ...\wedge dt_{q_X}^{(j)}$, 
and the integral becomes  
$$
\sum\limits_j 
\int\limits_{ \{ w\in A_z(X)|\tau(w,0)\le \delta \} \cap U_j }
e^{ -(n+1)k\ln \cosh \frac{\tau (w,0)}{2} }\psi^{(j)}(t)f^{(j)}(t)  dt_1^{(j)}\wedge ...\wedge dt_{q_X}^{(j)}.
$$ 
Now we will work with the integral 
\begin{equation}
\label{intlaplace}
\int\limits_{ \{ w\in A_z(X)|\tau(w,0)\le \delta \} \cap U_j }
e^{ -\frac{(n+1)k}{2}\ln \cosh^2 \frac{\tau (w,0)}{2} }
\psi^{(j)}(t)f^{(j)}(t)  
dt_1^{(j)}...dt_{q_X}^{(j)} 
\end{equation}
Apply the multivariable Laplace method.   
If the point $w=0$ is in $U_j$ or on the boundary of $U_j$, then 
the appropriate statement is, respectively, 
Theorem 3 p. 495  or (5.15) p. 498 in \cite{wong:89}. 
If the point  $w=0$ is not in $\overline{U_j}$,
then it follows that the contribution from the $j$-th 
integral is rapidly decreasing as $k\to \infty$,
by an argument similar to the one that has already 
been used earlier. 

In order to use the Laplace method, 
we need to show that the Hessian matrix $H_z$ of the function  
$\ln \cosh ^2 \frac{\tau (w,0)}{2}=-\ln (-\langle w,w\rangle )$  at $w=0$ 
is positive definite. We have: for $l\in \{ 1,...,q_X\}$, $p\in \{1,...,q_X\}$
$$
\frac{\partial}{\partial t_p}( -\ln (-\langle w,w\rangle )=\frac{1}{-\langle w,w\rangle }
\sum\limits_{r=1}^n \Bigl ( w_r \frac{\partial \bar{w}_r}{\partial t_p}+\bar{w}_r \frac{\partial w_r}{\partial t_p}\Bigr )
$$
$$
\frac{\partial ^2}{\partial t_l\partial t_p}( -\ln (-\langle w,w\rangle )\Bigl | _{w=0}
=\sum\limits_{r=1}^n \Bigl ( \frac{\partial w_r}{\partial t_l} \frac{\partial \bar{w}_r}{\partial t_p}+
\frac{\partial \bar{w}_r}{\partial t_l} \frac{\partial w_r}{\partial t_p}\Bigr ).
$$
Therefore $H_{\zeta}=B_{\zeta}\bar{B}_{\zeta}^T+\bar{B}_{\zeta}B_{\zeta}^T$, where $B_{\zeta}$ is the $q_X\times n$ 
matrix $\begin{pmatrix} \frac{\partial w_1}{\partial t_1} & ... & \frac{\partial w_n}{\partial t_1} \\
... & & ... \\  
\frac{\partial w_1}{\partial t_{q_X}} & ... & \frac{\partial w_n}{\partial t_{q_X}}
\end{pmatrix}$.
The matrix 
$H_{\zeta}$ is symmetric. The matrices $H_{\zeta}$, $B_{\zeta}\bar{B}_{\zeta}^T$, $\bar{B}_{\zeta}B_{\zeta}^T$ 
are positive semidefinite, because for a vector $v\in \C^{q_X}$ 
$(B_{\zeta}\bar{B}_{\zeta}^Tv)^T\bar{v}=(\bar{B}_{\zeta}^T v)^T\overline{\bar{B}_{\zeta}^T v}$ and  
$(\bar{B}_{\zeta}B_{\zeta}^Tv)^T\bar{v}=(B_{\zeta}^T v)^T\overline{B_{\zeta}^T v}$. 
It remains to show: if $H_{\zeta}v=0$, then $v=0$. If  $H_{\zeta}v=0$, then
$\bar{v}^TH_{\zeta}v=0$, and it follows that 
$B_{\zeta}^T v=0$. But $rk (B_{\zeta}^T)=q_X$, hence $\dim \ker B_{\zeta}^T=0$, therefore $v=0$. 
Thus $H_{\zeta}$ is positive definite. 

If the point $w=0$ is in $U_j$, then the integral (\ref{intlaplace}) is asymptotic, as $k\to\infty$, to 
$\Bigl (\frac{4\pi}{(n+1)k}\Bigr )^{\frac{q_X}{2}}
 \psi^{(j)}f^{(j)} \Bigr| _{w=0}
(\det H_{\zeta})^{-\frac{1}{2}}$, 
and if $0$ is on the boundary of $U_j$, then the integral (\ref{intlaplace}) is asymptotic to 
$\frac{1}{2}\Bigl (\frac{4\pi}{(n+1)k}\Bigr )^{\frac{q_X}{2}}
 \psi^{(j)}f^{(j)}\Bigr| _{w=0}
(\det H_{\zeta})^{-\frac{1}{2}}$.  
We conclude: 
$$
c(\B^n ,k)
\int\limits_{Y}\int\limits_{ \{ \zeta\in X|\tau(z,\zeta)\le \delta \}  }
e^{ -(n+1)k\ln \cosh \frac{\tau (z,\zeta)}{2} }\nu_{X}(\zeta)\nu_{Y}(z)=
$$
$$
c(\B^n ,k)\int\limits_{Y} \sum\limits_j 
\int\limits_{ \{ w\in A_z(X)|\tau(w,0)\le \delta \} \cap U_j }
e^{ -\frac{(n+1)k}{2}\ln \cosh^2 \frac{\tau (w,0)}{2}  }
\psi^{(j)}(t)f^{(j)}(t)  dt_1^{(j)}... dt_{q_X}^{(j)} \ \nu_{Y}(z)\sim 
$$
$$
c(\B^n ,k)Ck^{-\frac{q_X}{2}}
$$
and the statements $(i)$, $(ii)$ now follow 
from Remark \ref{remconstball}. 
For the constant $C$ we have: $C>0$, because the number
$f^{(j)}\Bigr| _{w=0}(\det H_{\zeta})^{-\frac{1}{2}}$ is positive for each $j$,
the value of the function $\psi^{(j)}$ at the point $w=0$ is nonnegative for
each $j$, and there is $j_0$ such that the point $w=0$ is in $U_{j_0}$ and
$ \psi^{(j_0)}\Bigr| _{w=0}>0$.
$\Box$
\begin{remark}
The remainder in Theorem \ref{thxx}$(ii)$ is determined by $I_2$, $I_1^{(2)}$, the error term in the Laplace approximation 
and the error in the Stirling formula. 
\end{remark}
In the examples below, for specific $X$ and $Y$, we will work out the integral   
$$
I_1=
c(\B^n ,k)  
\int_Y\int_X 
\frac{(\langle z,z\rangle\langle \zeta,\zeta\rangle) ^{ \frac{(n+1)k}{2}  }}
{(-\langle z,\zeta\rangle )  
^{ (n+1)k } }\nu_{X}(\zeta)\nu_{Y}(z).  
$$
This term appeared in the proof of Theorem \ref{thxx} as the term that determines 
the behaviour of 
$(\Theta_X^{(j;k)},\Theta_Y^{(j;k)})$ as $k\to\infty$. 

\begin{example} 
\label{exlinesegment}
Let $Y=X\subset \B^n$ be a ($1$-dimensional) line segment defined by 
$z_1=te^{i\varphi}$, $-\alpha < t < \alpha$, 
where $\alpha\in (0,1)$ and $\varphi\in [0,\frac{\pi}{2}]$ are fixed, $z_j=0$ for $j>1$, and let $\nu_X=dt$. 
If $n=1$, then $X$ is a Lagrangian submanifold of $\B^1$. For arbitrary $n$ such $X$ is totally real.  
We have: 
$$
I_1=
c(\B^n ,k)  
\int\limits_{-\alpha}^{\alpha}
\int\limits_{-\alpha}^{\alpha}
\Bigl ( \frac{(1-t^2)(1-T^2) }{ (1-tT)^2 } 
\Bigr )^{ \frac{(n+1)k}{2}  } dt \ dT. 
$$
Here $\zeta=Te^{i\varphi}$. For a fixed $T$ 
denote $f(t)=\Bigl ( \frac{(1-t^2)(1-T^2) }{ (1-tT)^2 }\Bigr ) ^{1/2}$. We have: 
$$
f(T)=1, \  \frac{df}{dt}\Bigr | _{t=T}=0, \    \frac{d^2f}{dt^2}\Bigr | _{t=T}=-\frac{1}{(1-T^2)^2}<0. 
$$
The function $f$ has a maximum at $t=T$. 
Applying the $1$-dimensional Laplace approximation  ((1.5) p. 60 \cite{wong:89} or (5.1.21) \cite{bleistein:75}) 
we get: as $k\to\infty$ 
$$
\int\limits_{-\alpha}^{\alpha} f(t)^{ (n+1)k } dt\sim 
\Bigl ( \frac{-2\pi}{ (n+1)kf''(T)  }  \Bigr ) ^{ \frac{1}{2}}=
\sqrt{\frac{2\pi}{(n+1)k}}(1-T^2), 
$$
hence 
$$
I_1\sim c(\B^n ,k) \sqrt{\frac{2\pi}{(n+1)k}}\int\limits_{-\alpha}^{\alpha}(1-T^2)dT
=c(\B^n ,k) \sqrt{\frac{2\pi}{(n+1)k}}2(\alpha-\frac{\alpha^3}{3})\sim 
c(n)(\alpha-\frac{\alpha^3}{3}) k^{n-\frac{1}{2}},
$$
where $c(n)=\frac{(n+1)^{n-\frac{1}{2} } }{n!}2\sqrt{2\pi}$. 
\end{example}
\begin{example} 
\label{excircle}
Let $Y=X\subset \B^n$ ($n\ge 2$) be the circle of radius $0<\alpha <1$ in the $x_1x_2$-plane
centered at $(x_1,x_2)=(0,0)$: 
$z_1=x_1=\alpha\cos \Theta$, $z_2=x_2=\alpha\sin \Theta$, $0\le \Theta <2\pi$, 
$y_1=y_2=0$, $z_j=0$ for $j>2$.
Let $\nu_X=d\Theta$. For arbitrary $n\ge 2$ such $X$ is totally real.  
We have: 
$$
I_1=c(\B^n ,k)\int_X\int_X
\frac{[(1-x_1^2-x_2^2)(1-(Re(\zeta_1))^2-(Re(\zeta_2))^2)]^{\frac{(n+1)k}{2}}}
     {(1-x_1Re(\zeta_1)-x_2Re(\zeta_2))^{(n+1)k}}
 d\Theta \ d\varphi=
$$
$$
c(\B^n ,k)  
\int\limits_{0}^{2\pi}
\int\limits_0^{2\pi }
\Bigl ( \frac{1-\alpha^2 }{ 1-\alpha^2\cos(\Theta-\varphi)} 
\Bigr )^{ (n+1)k  } d\Theta \ d\varphi. 
$$
Here $\zeta_1=\alpha\cos \varphi$, $0\le \varphi <2\pi$, 
$\zeta_2=\alpha\sin \varphi$, $Im(\zeta_1)=Im(\zeta_2)=0$, $\zeta_j=0$ for $j>2$.  
For a fixed $\varphi$ 
denote $f(\Theta)= \frac{1-\alpha^2 }{ 1-\alpha^2\cos(\Theta-\varphi)} $. We have: 
$$
f(\varphi)=1, \  \frac{df}{d\Theta}\Bigr | _{\Theta=\varphi}=0, \    
\frac{d^2f}{d\Theta^2}\Bigr | _{\Theta=\varphi}=-\frac{\alpha^2}{1-\alpha^2}<0,
$$
$f$ has a local maximum at $\Theta=\varphi$. 
Applying the $1$-dimensional Laplace approximation ((1.5) p. 60 \cite{wong:89} or (5.1.21) \cite{bleistein:75}) 
we get:
$$
\int\limits_0^{2\pi} f(\Theta)^{ (n+1)k } d\Theta\sim 
\Bigl ( \frac{-2\pi}{ (n+1)kf''(\varphi)  }  \Bigr ) ^{ \frac{1}{2}}=
\sqrt{\frac{2\pi}{(n+1)k}}\frac{\sqrt{1-\alpha^2}}{\alpha}, 
$$
hence 
$$
I_1\sim c(\B^n ,k) \sqrt{\frac{2\pi}{(n+1)k}}\frac{\sqrt{1-\alpha^2}}{\alpha}
\int\limits_{0}^{2\pi}d\varphi
=c(\B^n ,k) \sqrt{\frac{2\pi}{(n+1)k}}\frac{\sqrt{1-\alpha^2}}{\alpha} 2\pi
\sim 
k^{n-\frac{1}{2}}c(n)\frac{\sqrt{1-\alpha^2}}{\alpha},
$$
where $c(n)=\frac{(n+1)^{n-\frac{1}{2} } }{n!}2\pi\sqrt{2\pi}$. 
\end{example}
\begin{example} 
\label{exdisc}  
Let $Y=X\subset \B^n$ ($n\ge 2$) be the disc of radius $\alpha\in (0,1)$ in the $x_1x_2$-plane
centered at $(x_1,x_2)=(0,0)$. Thus, $X$ is defined by 
$x_1^2+x_2^2<\alpha^2$,   $y_1=y_2=0$, 
$z_j=0$ for $j>2$. Let $\nu_X=dx_1\wedge dx_2$. 
For arbitrary $n\ge 2$ such $X$ is totally real.  If $n=2$, then $X$ is a Lagrangian submanifold of $\B^2$.
$$
I_1=
c(\B^n ,k)
\int\limits_{X}\int\limits_X
\Bigl ( \frac{(1-x_1^2-x_2^2)(1-u_1^2-u_2^2)}{(1-x_1u_1-x_2u_2)^2} \Bigr ) ^{ \frac{(n+1)k}{2}  }dx_1dx_2du_1du_2=
$$
$$
c(\B^n ,k)
\int\limits_{X}\int\limits_X e^{\frac{(n+1)k}{2}\ln
\frac{(1-x_1^2-x_2^2)(1-u_1^2-u_2^2)}{(1-x_1u_1-x_2u_2)^2} }dx_1dx_2du_1du_2,
$$
where $u_1=Re(\zeta_1)$, $u_2=Re(\zeta_2)$. For fixed $u_1$, $u_2$ let 
$f(x_1, x_2)=- \ln \frac{(1-x_1^2-x_2^2)(1-u_1^2-u_2^2)}{(1-x_1u_1-x_2u_2)^2}$. 
We have: $f(u_1,u_2)=0$,
$$
\frac{\partial f}{\partial x_j}=2\Bigl ( \frac{x_j}{1-x_1^2-x_2^2}-\frac{u_j}{1-x_1u_1-x_2u_2}\Bigr ) , \ j=1,2 
$$
$$
\frac{\partial f}{\partial x_1} \Bigr |_{(u_1,u_2)}=\frac{\partial f}{\partial x_2} \Bigr |_{(u_1,u_2)}=0
$$
$$
\frac{\partial ^2 f}{\partial x_j^2}=2\Bigl ( \frac{1-x_1^2-x_2^2+2x_j^2}{(1-x_1^2-x_2^2)^2}-
\frac{u_j^2}{(1-x_1u_1-x_2u_2)^2} \Bigr ) , \ j=1,2
$$
$$
\frac{\partial ^2 f}{\partial x_1 \partial x_2}=2\Bigl ( \frac{2x_1x_2}{(1-x_1^2-x_2^2)^2}-
\frac{u_1u_2}{(1-x_1u_1-x_2u_2)^2} \Bigr )
$$
$$
\frac{\partial ^2 f}{\partial x_1^2} \Bigr |_{(u_1,u_2)}=2\frac{1-u_2^2}{(1-u_1^2-u_2^2)^2}>0
$$
$$
H(u_1,u_2)= \Bigl ( 
\frac{\partial ^2 f}{\partial x_1^2} \frac{\partial ^2 f}{\partial x_2^2} 
-\Bigl ( \frac{\partial ^2 f}{\partial x_1 \partial x_2} \Bigr )^2\Bigr )\Bigr |_{(u_1,u_2)}=
\frac{4}{(1-u_1^2-u_2^2)^3}>0
$$
Using Laplace approximation in $\R^2$ (\cite{wong:89} p. 495 or Theorem 2 \cite{hsu:51}) we get: 
for a fixed $\zeta$ 
$$
\int\limits_X e^{-\frac{(n+1)k}{2}f(x_1,x_2) }
dx_1dx_2 \sim 
\frac{4\pi}{(n+1)k\sqrt{H(u_1,u_2)}}=\frac{2\pi}{(n+1)k}(1-u_1^2-u_2^2)^{\frac{3}{2}}
$$
and
$$
I_1\sim 
c(\B^n ,k)\frac{2\pi}{(n+1)k}\int\limits_{X}(1-u_1^2-u_2^2)^{\frac{3}{2}}du_1du_2=
c(\B^n ,k)\frac{4\pi ^2}{5}\frac{1}{(n+1)k}(1-(1-\alpha^2)^{\frac{5}{2}})\sim 
$$
$$
k^{n-1}\frac{4\pi ^2}{5} \frac{(n+1)^{n-1}}{n!}(1-(1-\alpha^2)^{\frac{5}{2}}). 
$$
\end{example}
\begin{example} 
  \label{exlinedisc}
  Let $\beta\in (0,1)$ and $\alpha\in (\beta,1)$ be fixed. 
  Let $Y$ be the line segment in $\B^n$ ($n\ge 2$) defined by $-\beta <x_1<\beta$,
  $y_1=0$, $z_j=0$ for $j>1$, and let 
$X$ be the disc defined by $x_1^2+x_2^2<\alpha ^2$,  $y_1=y_2=0$, 
$z_j=0$ for $j>2$. Let $\nu_X=dx_1\wedge dx_2$ and $\nu_Y=dx_1$. 
$$
I_1=
c(\B^n ,k)
\int\limits_{Y}\int\limits_X
\Bigl ( \frac{(1-x_1^2-x_2^2)(1-u_1^2)}{(1-x_1u_1)^2} \Bigr ) ^{ \frac{(n+1)k}{2}  }dx_1dx_2du_1=
$$
$$
c(\B^n ,k)
\int\limits_{Y}\int\limits_X e^{\frac{(n+1)k}{2}\ln
\frac{(1-x_1^2-x_2^2)(1-u_1^2)}{(1-x_1u_1)^2} }dx_1dx_2du_1,
$$
where $u_1=Re(\zeta_1)$. For a fixed $u_1$ let 
$f(x_1, x_2)=- \ln \frac{(1-x_1^2-x_2^2)(1-u_1^2)}{(1-x_1u_1)^2}$. 
We have: $f(u_1,0)=0$,
$$
\frac{\partial f}{\partial x_1}=2\Bigl ( \frac{x_1}{1-x_1^2-x_2^2}-\frac{u_1}{1-x_1u_1}\Bigr ) , \  
\frac{\partial f}{\partial x_2}=\frac{2x_2}{1-x_1^2-x_2^2}
$$
$$
\frac{\partial f}{\partial x_1} \Bigr |_{(u_1,0)}=\frac{\partial f}{\partial x_2} \Bigr |_{(u_1,0)}=0
$$
$$
\frac{\partial ^2 f}{\partial x_1^2}=2\Bigl ( \frac{1+x_1^2-x_2^2}{(1-x_1^2-x_2^2)^2}-
\frac{u_1^2}{(1-x_1u_1)^2} \Bigr ) , \ 
\frac{\partial ^2 f}{\partial x_2^2}=2\frac{1-x_1^2+x_2^2}{(1-x_1^2-x_2^2)^2}
$$
$$
\frac{\partial ^2 f}{\partial x_1 \partial x_2}=\frac{-4x_1x_2}{(1-x_1^2-x_2^2)^2}
$$
$$
\frac{\partial ^2 f}{\partial x_1^2} \Bigr |_{(u_1,0)}=\frac{2}{(1-u_1^2)^2}>0
$$
$$
H(u_1,0)= \Bigl ( 
\frac{\partial ^2 f}{\partial x_1^2} \frac{\partial ^2 f}{\partial x_2^2} 
-\Bigl ( \frac{\partial ^2 f}{\partial x_1 \partial x_2} \Bigr )^2\Bigr )\Bigr |_{(u_1,0)}=
\frac{4}{(1-u_1^2)^3}>0
$$
Using Laplace approximation in $\R^2$ (\cite{wong:89} p. 495 or Theorem 2 \cite{hsu:51}) we get: 
for a fixed $\zeta$ 
$$
\int\limits_X e^{-\frac{(n+1)k}{2}f(x_1,x_2) }
dx_1dx_2 \sim 
\frac{4\pi}{(n+1)k\sqrt{H(u_1,0)}}=\frac{2\pi}{(n+1)k}(1-u_1^2)^{\frac{3}{2}}
$$
and as $k\to\infty$ 
$$
I_1\sim 
c(\B^n ,k)\frac{2\pi}{(n+1)k}\int_{-\beta}^{\beta}(1-u_1^2)^{\frac{3}{2}}du_1
\sim {\mathrm{const}}(n,\beta ) \ k^{n-1}.
$$
\end{example}

\begin{example} 
\label{3ball}
\noindent Let $\alpha\in (0,1)$ be fixed. Define a submanifold of $\B^n$, $n\ge 2$, by  
$$
X=\{ (x_1,y_1,...,x_n,y_n)\in \B^n \ | \  x_1^2+y_1^2+x_2^2<\alpha^2, x_2> 0, y_2=0, x_j=y_j=0  \
{\mathrm{for}} \ j>2\}. 
$$
and set  $\nu_X=dx_1\wedge dy_1\wedge dx_2$.

As a remark, we point out that $X$ is a CR submanifold which is not totally real and not complex
(see subsection \ref{appsec3} of the Appendix).

We will now estimate $I_1$, with $X=Y$:
$$
I_1=c(\B^n ,k)  
\int_X \int_X 
\frac{(\langle z,z\rangle\langle \zeta,\zeta\rangle) ^{ \frac{(n+1)k}{2}  }}
{(-\langle z,\zeta\rangle )  
^{ (n+1)k } }\nu_{X}(z)\nu_{X}(\zeta).
$$
We set all coordinates, except for $x_1$, $y_1$ and $x_2$, to be zero. We 
will use the spherical coordinates:
$$
x_1=\rho \sin\Phi \cos \Theta
$$
$$
y_1=\rho \sin\Phi\sin\Theta
$$
$$
x_2=\rho \cos\Phi
$$
$$
0< \rho< \alpha, \ 0\le \Theta<2\pi, \ 0\le \Phi < \frac{\pi}{2}
$$
$$
u_1=r \sin\psi \cos \beta
$$
$$
v_1=r \sin\psi\sin\beta
$$
$$
u_2=r \cos\psi
$$
$$
0< r< \alpha, \ 0\le \beta<2\pi, \ 0\le \psi< \frac{\pi}{2},
$$
where $\zeta_1=u_1+iv_1$ ($u_1,v_1\in\R$), $Re(\zeta_2)=u_2$. 
The integral becomes 
$$
I_1=c(\B^n,k) \int_X\int_X 
\frac{[(1-x_1^2-y_1^2-x_2^2)(1-u_1^2-v_1^2-u_2^2)]^{\frac{(n+1)k}{2}}}{(1-(x_1+iy_1)(u_1-iv_1)-x_2u_2)^{(n+1)k}}
dx_1 dy_1 dx_2 du_1 dv_1 du_2.
$$
For fixed $u_1$, $v_1$, $u_2$, the integral  
$$
\int_X 
\frac{(1-x_1^2-y_1^2-x_2^2)^{\frac{(n+1)k}{2}}}{(1-(x_1+iy_1)(u_1-iv_1)-x_2u_2)^{(n+1)k}}
dx_1 dy_1 dx_2=
$$
$$
\int_X \frac{(1-\rho^2)^{\frac{(n+1)k}{2}} \rho^2\sin\Phi }
{(1-\rho\sin\Phi(\cos\Theta+i\sin\Theta)(u_1-iv_1)-u_2\rho\cos\Phi )^{(n+1)k}}
d\rho  \ d\Phi  \ d\Theta =
$$ 
$$
\int_0^{\frac{\pi}{2}}\sin\Phi \int_0^{\alpha} (1-\rho^2)^{\frac{(n+1)k}{2}} \rho^2 
\int_0^{2\pi} \frac{1}
    {(1-\rho\sin\Phi  \ e^{i\Theta}(u_1-iv_1)-u_2\rho\cos\Phi )^{(n+1)k}ie^{i\Theta}}
    $$
    $$
    d(e^{i\Theta})d\rho  \ d\Phi  =
2\pi \int_0^{\frac{\pi}{2}}\sin\Phi \int_0^{\alpha} (1-\rho^2)^{\frac{(n+1)k}{2}} \rho^2 
\frac{1}{(1-u_2\rho\cos\Phi )^{(n+1)k}} d\rho  \ d\Phi .
$$
Apply the Laplace method (\cite{wong:89} Theorem 3 p. 495)  to the integral 
$$
\int_0^{\alpha} (1-\rho^2)^{\frac{(n+1)k}{2}} \rho^2 
\frac{1}{(1-u_2\rho\cos\Phi)^{(n+1)k}} d\rho=
\int_0^{\alpha} e^{-(n+1)kf(\rho)} \rho^2  d\rho, 
$$
where
$$
f(\rho)=-\ln \frac{\sqrt{1-\rho^2}}
{1-u_2\rho\cos\Phi }.
$$
We have:
$$
\frac{df}{d\rho}=-\frac{\rho-u_2\cos \Phi}{(1-\rho^2)(1-u_2\rho\cos\Phi )}, 
$$
$$
\frac{d^2f}{d\rho^2}\Bigr |_{\rho=u_2\cos\Phi}=\frac{1}{(1-(u_2\cos\Phi )^2)^2}>0, 
$$
and as $k\to\infty$ 
$$
\int_0^{\alpha} e^{-(n+1)kf(\rho)} \rho^2  d\rho \sim
e^{-(n+1)kf(u_2\cos\Phi)}\Bigl [ 
\sqrt{\frac{2\pi}{(n+1)kf''(u_2\cos\Phi)}}
(u_2\cos\Phi)^2 +c_1((n+1)k)^{-\frac{3}{2}}+
$$
$$
c_2((n+1)k)^{-\frac{5}{2}}+...\Bigr ]
$$
where the constants $c_1$, $c_2$,..., depend on $\Phi$ and $u_2$. So, 
$$
\int_0^{\alpha} e^{-(n+1)kf(\rho)} \rho^2  d\rho \sim (1-(u_2\cos\Phi)^2)^{-\frac{(n+1)k}{2}}
\Bigl [ 
  \sqrt{\frac{2\pi}{(n+1)k}}(u_2\cos\Phi)^2(1-(u_2\cos\Phi)^2)+
  $$
  $$
   c_1((n+1)k)^{-\frac{3}{2}}+
c_2((n+1)k)^{-\frac{5}{2}}+...\Bigr ] .
$$  
Therefore
$$
I_1\sim c(\B^n,k)2\pi  \int_X \int_0^{\frac{\pi}{2}}
 (1-(u_2\cos\Phi)^2)^{-\frac{(n+1)k}{2}}
\Bigl [ 
  \sqrt{\frac{2\pi}{(n+1)k}}(u_2\cos\Phi)^2(1-(u_2\cos\Phi)^2)+
  $$
  $$
   c_1((n+1)k)^{-\frac{3}{2}}+
c_2((n+1)k)^{-\frac{5}{2}}+...\Bigr ] 
(1-u_1^2-v_1^2-u_2^2)^{\frac{(n+1)k}{2}}\sin\Phi 
 \ d\Phi  \ du_1 dv_1 du_2=
$$
$$
c(\B^n,k)2\pi \int_0^{2\pi}d\beta
\int_0^{\frac{\pi}{2}} \int_0^{\alpha}
\int_0^{\frac{\pi}{2}}
(1-(r\cos\psi\cos\Phi)^2)^{-\frac{(n+1)k}{2}}
$$
$$
\Bigl [ 
  \sqrt{\frac{2\pi}{(n+1)k}}(r\cos\psi\cos\Phi)^2(1-(r\cos\psi\cos\Phi)^2)+
   c_1((n+1)k)^{-\frac{3}{2}}+
   c_2((n+1)k)^{-\frac{5}{2}}+...\Bigr ]
$$
$$
(1-r^2)^{\frac{(n+1)k}{2}}\sin\Phi \  
 r^2\sin\psi \ 
d\Phi  \ dr \  d\psi=
$$
$$
c(\B^n,k)4\pi^2 
\int_0^{\frac{\pi}{2}}
\int_0^{\frac{\pi}{2}} \int_0^{\alpha}
(1-(r\cos\psi\cos\Phi)^2)^{-\frac{(n+1)k}{2}}
$$
$$
\Bigl [ 
  \sqrt{\frac{2\pi}{(n+1)k}}(r\cos\psi\cos\Phi)^2(1-(r\cos\psi\cos\Phi)^2)+
   c_1((n+1)k)^{-\frac{3}{2}}+
   c_2((n+1)k)^{-\frac{5}{2}}+...\Bigr ]
$$
$$
(1-r^2)^{\frac{(n+1)k}{2}}\sin\Phi \  
 r^2\sin\psi \ 
dr  \ d\Phi \   d\psi =
$$
$$
c(\B^n,k)4\pi^2  \sqrt{\frac{2\pi}{(n+1)k}}
\int_0^{\frac{\pi}{2}} (\cos\psi )^2\sin\psi 
\int_0^{\frac{\pi}{2}} (\cos\Phi)^2\sin\Phi   \int_0^{\alpha}
(1-(r\cos\psi\cos\Phi)^2)^{-\frac{(n+1)k}{2}} 
$$
$$
(1-(r\cos\psi\cos\Phi)^2)
(1-r^2)^{\frac{(n+1)k}{2}}   
 r^4 
dr  \ d\Phi \  d\psi +
$$
$$
c(\B^n,k)4\pi^2 
\int_0^{\frac{\pi}{2}}\sin\psi 
\int_0^{\frac{\pi}{2}} \sin\Phi \int_0^{\alpha}
(1-(r\cos\psi\cos\Phi)^2)^{-\frac{(n+1)k}{2}}
$$
$$
\Bigl [ 
   c_1((n+1)k)^{-\frac{3}{2}}+
   c_2((n+1)k)^{-\frac{5}{2}}+...\Bigr ]
(1-r^2)^{\frac{(n+1)k}{2}}   
 r^2  
dr \ d\Phi \ d\psi .
$$
For fixed values of $\Phi$, $\Psi$, apply the Laplace method
 (\cite{wong:89} Theorem 3 p. 495) to the integral
$$
\int_0^{\alpha}e^{-(n+1)kf(r)}  
(1-(r\cos\psi\cos\Phi)^2)
r^4 dr=\frac{1}{2}
\int_{-\alpha}^{\alpha}e^{-(n+1)kf(r)}  
(1-(r\cos\psi\cos\Phi)^2)
r^4 dr
$$
 where
$$
f(r)=-\frac{1}{2}\ln\frac{1-r^2}{1-(r\cos\psi\cos\Phi)^2}.
$$
We get:
$$
\frac{df}{dr}=\frac{r-r(\cos\Phi\cos\psi)^2}{(1-r^2)(1-(r\cos\Phi\cos\psi)^2)}
$$
$$
\frac{d^2f}{dr^2}\Bigr |_{r=0}=1-(\cos\Phi\cos\psi)^2)>0 
$$
and as $k\to\infty$ 
$$
\int_{-\alpha}^{\alpha}e^{-(n+1)kf(r)}  
(1-(r\cos\psi\cos\Phi)^2)
r^4 dr\sim
$$
$$
e^{-(n+1)kf(0)}\Bigl [ \sqrt{\frac{2\pi}{(n+1)kf''(0)}}
(1-(r\cos\psi\cos\Phi)^2)r^4\Bigr |_{r=0}+O(k^{-\frac{3}{2}})\Bigr ].
$$
We also have, for fixed values of $\Phi$, $\psi$:
$$
| \int_0^{\alpha}
(1-(r\cos\psi\cos\Phi)^2)^{-\frac{(n+1)k}{2}}
\Bigl [ 
   c_1((n+1)k)^{-\frac{3}{2}}+
   c_2((n+1)k)^{-\frac{5}{2}}+...\Bigr ]
(1-r^2)^{\frac{(n+1)k}{2}}
r^2dr|
$$
$$
\le {\mathrm{const}} \ k^{-3},
$$
since by the Laplace method, as above,  
$$
\int_0^{\alpha}
(1-(r\cos\psi\cos\Phi)^2)^{-\frac{(n+1)k}{2}}
(1-r^2)^{\frac{(n+1)k}{2}}
r^2dr=\frac{1}{2}
\int_{-\alpha}^{\alpha}
e^{-(n+1)k\frac{1}{2}\ln \frac{1-(r\cos\psi\cos\Phi)^2}{1-r^2}}
r^2dr\sim
$$
$$
O(k^{-\frac{3}{2}}).
$$
We recall: $c(\B^n ,k)\sim \frac{(n+1)^n}{n!}k^n\Bigl (1+O(\frac{1}{k})\Bigr )$
(Remark \ref{remconstball}). Combining all this together, we conclude:
$$
|I_1|\le {\mathrm{const}} \ k^{n-2}
$$
It follows that
$$
(\Theta_X^{(j;k)},\Theta_X^{(j;k)})\le {\mathrm{const}} \ k^{n-2}.
$$
We have: ${\mathrm{dim}}_{\R}X=3$,
for $n>2$ we are in the setting of  part {\it (i)} of Theorem \ref{thxx}, 
the asymptotic inequality holds and it is a strict inequality, 
and the asymptotic behaviour we observe here is different 
from the asymptotics for totally real submanifolds in
Theorem \ref{thxx}{\it (ii)}, for which
$(\Theta_X^{(j;k)},\Theta_X^{(j;k)})$ would have been asymptotic
to $Ck^{n-\frac{{\mathrm{dim}}_{\R}X}{2}}$ with $C>0$.
\end{example}
\begin{remark}
  In  \cite{borthwick:95} it is shown that 
  the square of the norm of the sections of powers of
  the line bundle associated to  a  Bohr-Sommerfeld Lagrangian submanifold
of a compact $n$-dimensional K\"ahler manifold
 grows as a constant times 
 $k^{\frac{n}{2}}$. In \cite{ioos:18} an analogous result for
 a $q$-dimensional isotropic Bohr-Sommerfeld submanifold of a $2n$-dimensional
 symplectic manifold and associated sections of vector bundles gives 
 the leading term of $const \ k^{n-\frac{q}{2}}$. Our 
 Theorem \ref{thxx}{\it (ii})  gives the same leading term for norms of
 vector-valued automorphic forms associated to isotropic submanifolds of ball quotients. 
In the Example \ref{3ball} we have a submanifold $X$
which is not isotropic
and $(\Theta_X^{(j;k)},\Theta_X^{(j;k)})$ does not have the same kind of asymptotics. 
All this raises a general question
how the geometric properties of submanifolds are reflected in asymptotics
of the associated sections of vector bundles. 
\end{remark}

\appendix{}
\section{}
\subsection{Proof of Lemma \ref{propconvscalar}.} First, we observe that this statement is contained 
in the general framework of Chapters 5 and 7 of \cite{kollar:95}.  
Now we will present an actual proof. This is a modification of the proof of Prop. 1, p. 44 \cite{baily:73}, 
which will use that for a fixed $p\in D$  
$\int\limits _D|K(z,p)^2|dV_e(z)<\infty$. 
Let $A$ be a nonempty compact subset of $D$. We will prove that the series 
$\sum\limits_{\gamma\in\Gamma}\Bigl ( K(\gamma z,p) J(\gamma,z)\Bigr ) ^k$ 
converges absolutely and uniformly on $A$, for $k\ge2$. 
There are a compact subset $B$ of $D$ and $\delta>0$ such that $A$ is contained in the interior of $B$, 
and for any $a\in A$ there is a polydisc $P_a$ of Euclidean volume $\delta$, with center $a$, such that 
$\bar{P}_a\subset B$. Let $m_0$ be the number of elements in $\{ g\in \Gamma | gB\cap B\ne \emptyset\}$. 
First consider the case $k=2$. For $a\in A$ 
$$
|K(\gamma a,p)^2 J(\gamma, a)^2|\le \frac{1}{\delta}\int\limits_{P_a}|K(\gamma z,p)^2J(\gamma,z)^2|dV_e(z)=
\frac{1}{\delta}\int\limits_{\gamma P_a}|K(w,p)^2|dV_e(w)
$$
where $w=\gamma z$. We get: 
$$
\sum\limits_{\gamma\in\Gamma}|K(\gamma a,p)^2J(\gamma, a)^2|\le  
\frac{1}{\delta} \sum\limits_{\gamma\in\Gamma} \int\limits_{\gamma P_a}|K(w,p)|^2dV_e(w)\le
\frac{m_0}{\delta}\int\limits_D|K(w,p)|^2dV_e(w).
$$
The last inequality is justified by observing that if $\gamma P_a \cap \gamma' P_a \ne \emptyset$ 
for $\gamma,\gamma'\in \Gamma$,  
then $\gamma ^{-1}\gamma '\in \{ g\in\Gamma| gB \cap B\ne \emptyset\}$, so each $w\in D$ is in at most 
$m_0$ of the sets $\gamma P_a$, $\gamma\in \Gamma$. This settles the case $k=2$. 
Therefore for $a\in A$ $|K(\gamma a,p)J(\gamma ,a)|<1$ for all but at most finitely many $\gamma \in \Gamma$.  
When $|K(\gamma a,p)J(\gamma ,a)|<1$, $|K(\gamma a,p)J(\gamma ,a)|^k$ is a decreasing function of $k\ge 2$. 
The desired statement follows. 
$\Box$

\subsection{Proof of Lemma \ref{lemconstcdk}} 
\noindent From (\ref{constcdk}) for $D=\B ^n$ with $z=0$ we get: 
$$
c(\B^n ,k)\frac{n!}{\pi^n}\int\limits_{\B^n }(-\langle w,w\rangle )^{(n+1)(k-1)}(\frac{i}{2})^n 
dw_1\wedge d\bar{w}_1\wedge ...\wedge dw_n\wedge d\bar{w}_n=1. 
$$
The integral in the left hand side is equal to $\pi^n\frac{((n+1)(k-1))!}{((n+1)(k-1)+n)!}$, and the statement readily 
follows. To calculate this integral apply a change of variables 
$(w_1,\bar{w}_1)\to (R_1,\varphi_1)$, where $0\le R_1\le 1$, $0\le \varphi_1<2\pi$, 
$w_1=R_1e^{i\varphi_1}\sqrt{1-|w_2|^2-...-|w_n|^2}$. We get:
$$ 
\int\limits_{\B^n }(1-|w_1|^2-...-|w_n|^2 )^{(n+1)(k-1)}(\frac{i}{2})^n 
dw_1\wedge d\bar{w}_1\wedge ...\wedge dw_n\wedge d\bar{w}_n=
$$
$$
(\frac{i}{2})^{n-1} \int\limits_{\B^n }(1-|w_2|^2-...-|w_n|^2 )^{(n+1)(k-1)+1}(1-R_1^2)^{(n+1)(k-1)}R_1
dR_1 \wedge d\varphi_1\wedge dw_2\wedge d\bar{w}_2\wedge ...\wedge dw_n\wedge d\bar{w}_n, 
$$
then apply the change of variables 
$(w_2,\bar{w}_2)\to (R_2,\varphi_2)$, where $0\le R_2\le 1$, $0\le \varphi_2<2\pi$, 
$w_2=R_2e^{i\varphi_2}\sqrt{1-|w_3|^2-...-|w_n|^2}$ to transform the integral into 
$$
(\frac{i}{2})^{n-2} \int\limits_{\B^n }(1-|w_3|^2-...-|w_n|^2 )^{(n+1)(k-1)+2}(1-R_1^2)^{(n+1)(k-1)}R_1
(1-R_2^2)^{(n+1)(k-1)+1}R_2
dR_1\wedge d\varphi_1\wedge 
$$
$$
dR_2 \wedge d\varphi_2\wedge
dw_3\wedge d\bar{w}_3\wedge ...\wedge dw_n\wedge d\bar{w}_n, 
$$
and so on. At the end we get: 
$$ 
(2\pi)^{n-1}\int\limits_0^1(1-R_1^2)^{(n+1)(k-1)}R_1dR_1\int\limits_0^1(1-R_2^2)^{(n+1)(k-1)+1}R_2dR_2
... 
$$
$$
\int\limits_0^1(1-R_{n-1}^2)^{(n+1)(k-1)+{n-2}}R_{n-1}dR_{n-1}\int\limits_{|w_n|\le 1} (1-|w_n|^2)^{(n+1)(k-1)+(n-1)}
\frac{i}{2}dw_n\wedge d\bar{w}_n, 
$$
and with $w_n=R_ne^{i\varphi_n}$, $0\le R_n\le 1$, $0\le \varphi_n<2\pi$, the last integral is 
$$
2\pi\int\limits_0^1(1-R_n^2)^{(n+1)(k-1)+n-1}R_ndR_n.
$$ 
An elementary calculation now yields the answer. 
$\Box$

\subsection{Note for Example \ref{3ball}}
\label{appsec3}
In this subsection we explain why the submanifold 
$$
X=\{ (x_1,y_1,...,x_n,y_n)\in \B^n \ | \  x_1^2+y_1^2+x_2^2<\alpha^2, x_2> 0, y_2=0, x_j=y_j=0  \
{\mathrm{for}} \ j>2\}. 
$$
of $\B^n$, where  $n\ge 2$ and  $\alpha\in (0,1)$, is   
a CR submanifold, which is not totally real and not complex.
The relevant part of the complex hyperbolic metric (3.3)\cite{goldman:99} is,
up to a positive factor, 
$$
\frac{1}{(1-x_1^2-y_1^2-x_2^2-y_2^2)^2}[
  (1-x_2^2-y_2^2)(dx_1^2+dy_1^2)+(1-x_1^2-y_1^2)(dx_2^2+dy_2^2)
  $$
  $$
  +2(x_1x_2+y_1y_2)(dx_1dx_2+dy_1dy_2)+2(-y_1x_2+x_1y_2)(dx_1dy_2-dy_1dx_2)].
$$
The complex structure $J:T\B^n\to T\B^n$ acts as follows:
$$
\frac{\partial}{\partial x_j}\mapsto -\frac{\partial}{\partial y_j}, \
\frac{\partial}{\partial y_j}\mapsto \frac{\partial}{\partial x_j}
$$
for $j=1,2$. 

The distribution ${\mathcal{D}}=span\{ \frac{\partial}{\partial x_1} ,
\frac{\partial }{\partial y_1}\}$ is a holomorphic distribution on $X$. The
complementary orthogonal distribution 
${\mathcal{D}}^{\perp}=span\{ \frac{\partial}{\partial x_2}-
\frac{x_1x_2}{1-x_2^2}\frac{\partial}{\partial x_1}
-\frac{x_2y_1}{1-x_2^2}\frac{\partial}{\partial y_1}\}$
is a totally real distribution.
Indeed, $J( \frac{\partial}{\partial x_2}-
\frac{x_1x_2}{1-x_2^2}\frac{\partial}{\partial x_1}
-\frac{x_2y_1}{1-x_2^2}\frac{\partial}{\partial y_1})$ is orthogonal
to $ \dfrac{\partial}{\partial x_1}$,  $\dfrac{\partial}{\partial y_1}$,
$\dfrac{\partial}{\partial x_2}$, and is, therefore, in the normal bundle of $X$.


\begin{thebibliography}{999}

\bibitem{ali:07}S. Ali, M. Englis. 
\newblock {\em Matrix-valued Berezin-Toeplitz quantization.}
\newblock J. Math. Phys. 48 (2007), no. 5, 053504, 14pp.

\bibitem{alluhaibi:17}N. Alluhaibi.
  \newblock {\em On vector-valued automorphic forms on bounded symmetric domains.} 
    \newblock Ph.D. Thesis, University of Western Ontario, 2017. 

\bibitem{baily:73}W. Baily.  
\newblock {\em Introductory lectures on automorphic forms.} 
\newblock Iwanami Shoten, Publishers, Tokyo; Princeton University Press, Princeton, N.J., 1973. 

\bibitem{barron:12}T. Barron. 
\newblock {\em Quantization and automorphic forms.}  
\newblock Contemp. Math. 583 (2012), 211-219. 

\bibitem{barron:18}T. Barron.
\newblock {\em Closed geodesics and pluricanonical sections on ball quotients.} 
\newblock https://arxiv.org/abs/1808.01245

\bibitem{bell:67}D. Bell.
\newblock {\em Poincar\'e series representations of automorphic forms.}
\newblock Ph.D. Thesis, Brown University, 1967. 

\bibitem{bleistein:75}N. Bleistein, R. Handelsman. 
\newblock {\em Asymptotic expansions of integrals.}
\newblock Holt, Rinehart, Winston, 1975.

\bibitem{borel:65}A. Borel. 
\newblock {\em Introduction to automorphic forms.} In 
\newblock {\em Algebraic Groups and Discontinuous Subgroups (Proc. Sympos. Pure Math., Boulder, Colo., 1965) }
\newblock p. 199-210 AMS, Providence R.I. 1966. 

\bibitem{borthwick:95}D. Borthwick, T. Paul, A. Uribe. 
\newblock {\em Legendrian distributions with applications to relative Poincar\'e series.} 
\newblock Invent. Math. 122 (1995), no. 2, 359-402. 

\bibitem{burns:10}D. Burns, V. Guillemin, Z. Wang. 
\newblock {\em Stability functions.} 
\newblock Geom. Funct. Anal. 19 (2010), no. 5, 1258-1295. 

\bibitem{clery:13}F. Cl\'ery, G. van der Geer. 
\newblock {\em Generators for modules of vector-valued Picard modular forms.} 
\newblock Nagoya Math. J. 212 (2013), 19-57. 

\bibitem{debru:61}N. de Bruijn. 
\newblock {\em Asymptotic methods in analysis.}
\newblock Second edition. Bibliotheca Mathematica, Vol. IV. 
\newblock North-Holland Publishing Co., Amsterdam; P. Noordhoff Ltd., Groningen 1961. 

\bibitem{deber:06}M. Debernardi, R. Paoletti. 
\newblock {\em Equivariant asymptotics for Bohr-Sommerfeld Lagrangian submanifolds.}
\newblock Comm. Math. Phys. 267 (2006), no. 1, 227-263. 

\bibitem{foth:07}T. Foth. 
\newblock {\em Poincar\'e series on bounded symmetric domains.} 
\newblock Proc. AMS 135 (2007), no. 10, 3301-3308.

\bibitem{foth:08}T. Foth. 
\newblock {\em Legendrian tori and the semi-classical limit.} 
\newblock Diff. Geom. Appl. 26(2008), no.1, 63-74.

\bibitem{foth:01}T. Foth, S. Katok.
\newblock {\em Spanning sets for automorphic forms and dynamics of the frame flow on complex hyperbolic spaces.}
\newblock Ergodic Theory Dynam. Systems 21 (2001), no. 4, 1071-1099.

\bibitem{foth:04}T. Foth, S. Katok.
\newblock Appendix to {\em S. Katok.
Livshitz theorem for the unitary frame flow.} 
\newblock Ergodic Theory Dynam. Systems 24 (2004), no. 1, 127-140; pp. 137-140. 

\bibitem{freitag:14}E. Freitag, R. Manni.  
\newblock {\em Vector valued modular forms on three dimensional ball.} 
\newblock http://lanl.arxiv.org/abs/1404.3057

\bibitem{goldman:99}W. Goldman. 
\newblock {\em Complex hyperbolic geometry.}
\newblock Oxford Mathematical Monographs. Oxford Science Publications. 
\newblock The Clarendon Press, Oxford University Press, New York, 1999.

\bibitem{gorod:01}A. Gorodentsev, A. Tyurin. 
\newblock {\em Abelian Lagrangian algebraic geometry.}
\newblock   Izv. Math. 65 (2001), no. 3, 437-467.  

\bibitem{guillemin:15}V. Guillemin, A. Uribe, Z. Wang. 
\newblock {\em Semiclassical states associated with isotropic submanifolds of phase space.} 
\newblock   Lett. Math. Phys. 106 (2016), no. 12, 1695-1728. 

\bibitem{hsu:51}L. Hsu. 
\newblock {\em On the asymptotic evaluation of a class of multiple integrals involving a parameter.} 
\newblock Amer. J. Math. 73 (1951), 625-634. 

\bibitem{ioos:18}L. Ioos. 
  \newblock {\em   Quantization and isotropic submanifolds.}
  \newblock https://arxiv.org/abs/1802.09930

\bibitem{jeffrey:92}L. Jeffrey, J. Weitsman. 
\newblock {\em Bohr-Sommerfeld orbits in the moduli space of flat connections and the Verlinde dimension formula.} 
\newblock Comm. Math. Phys. 150 (1992), no. 3, 593-630. 

\bibitem{kato:84}S. Kato.
\newblock {\em A dimension formula for a certain space of automorphic forms of SU(p,1).} 
\newblock Math. Ann. 266 (1984), no. 4, 457-477. 

\bibitem{katok:85}S. Katok. 
\newblock {\em Closed geodesics, periods and arithmetic of modular forms.} 
\newblock Invent. Math. 80 (1985), no. 3, 469-480. 

\bibitem{katok:87}S. Katok, J. Millson. 
\newblock {\em Eichler-Shimura homology, intersection numbers and rational structures on spaces of modular forms.} 
\newblock Trans. Amer. Math. Soc. 300 (1987), no. 2, 737-757. 

\bibitem{kobayashi:87}S. Kobayashi. 
\newblock {\em Differential geometry of complex vector bundles.} 
\newblock Princeton University Press, Princeton, NJ; Princeton University Press, Princeton, NJ, 1987.

\bibitem{kojima:97}H. Kojima. 
\newblock {\em The formula for the dimension of the spaces of vector-valued holomorphic automorphic forms 
on the unitary group SU(1,p).} 
\newblock Kyushu J. Math. 51 (1997), no. 1, 57-76. 

\bibitem{kollar:95}J. Koll\'ar. 
\newblock {\em Shafarevich maps and automorphic forms.} 
\newblock Princeton University Press, Princeton, NJ, 1995. 

\bibitem{kudla:79}S. Kudla, J. Millson. 
\newblock {\em Harmonic differentials and closed geodesics on a Riemann surface.}
\newblock Invent. Math. 54 (1979), no. 3, 193-211. 

\bibitem{lu:15}Z. Lu, S. Zelditch.  
\newblock {\em Szeg\"o kernels and Poincar\'e series.} 
\newblock  J. Anal. Math. 130 (2016), 167-184.

\bibitem{ma:15}X. Ma, G. Marinescu. 
\newblock {\em Exponential estimate for the asymptotics of Bergman kernels.} 
\newblock Math. Ann. 362 (2015), no. 3-4, 1327-1347. 

\bibitem{narasimhan:64}M. Narasimhan, C. Seshadri.
  \newblock {\em  Holomorphic vector bundles on a compact Riemann surface.}
\newblock   Math. Ann. 155 (1964), 69-80.

\bibitem{paoletti:08}R. Paoletti. 
\newblock {\em A note on scaling asymptotics for Bohr-Sommerfeld Lagrangian submanifolds.}
\newblock Proc. Amer. Math. Soc. 136 (2008), no. 11, 4011-4017. 

\bibitem{parker:09}J. Parker. 
\newblock {\em Complex hyperbolic lattices.} 
\newblock In {\em Discrete groups and geometric structures,} 1-42,
\newblock Contemp. Math., 501, Amer. Math. Soc., Providence, RI, 2009. 

\bibitem{range:86}R. M. Range. 
\newblock {\em Holomorphic functions and integral representations in several complex variables.} 
\newblock Springer-Verlag, New York, 1986. 

\bibitem{rudin:80}W. Rudin. 
\newblock {\em Function theory in the unit ball of $\C^n$.}
\newblock  Springer-Verlag, New York-Berlin, 1980.

\bibitem{selberg:57}A. Selberg. 
\newblock {\em Automorphic functions and integral operators.} 
\newblock In {\em Collected papers}, vol. I,  Springer-Verlag, 1989; 464-468.

\bibitem{selberg:68}A. Selberg. 
\newblock {\em Recent developments in the theory of discontinuous groups of motions of symmetric spaces.} 
\newblock In {\em Collected papers}, vol. I,  Springer-Verlag, 1989; 546-567.

\bibitem{tong:83}Y. Tong, S. Wang. 
\newblock {\em Theta functions defined by geodesic cycles in quotients of SU(p,1).} 
\newblock Invent. Math. 71 (1983), no. 3, 467-499. 

\bibitem{wu:15}J. Wu, X. Wang. 
\newblock {\em  Poincar\'e series and very ampleness criterion for pluri-canonical bundles.}
  \newblock https://arxiv.org/abs/1504.00081
  
\bibitem{wong:89}R. Wong. 
\newblock {\em Asymptotic approximations of integrals.} 
\newblock Academic Press, Inc., Boston, MA, 1989. 

\bibitem{zhu:05}K. Zhu. 
\newblock {\em Spaces of holomorphic functions in the unit ball.} 
\newblock Springer-Verlag, New York, 2005. 


\end{thebibliography}
\end{document}